\documentclass{article}
\usepackage{amsmath}%
\usepackage{amsthm}%
\usepackage{amssymb}%

\usepackage{comment}

\newtheorem{theorem}{Theorem}[section]

\newtheorem{definition}[theorem]{Definition}
\newtheorem{proposition}[theorem]{Proposition}
\newtheorem{lemma}[theorem]{Lemma}
\newtheorem{remark}[theorem]{\it Remark}
\newtheorem{corollary}[theorem]{Corollary}
\newtheorem{example}[theorem]{\it Example}

\renewcommand{\indent}{\hspace*{5mm}}

\newcommand{\cF}{{\cal F}}
\newcommand{\cK}{{\cal K}}
\newcommand{\cP}{{\cal P}}
\newcommand{\cQ}{{\cal Q}}
\newcommand{\cR}{{\cal R}}
\newcommand{\cS}{{\cal S}}
\newcommand{\cT}{{\cal T}}

\newcommand{\fm}{{\mathfrak m}}

\renewcommand{\ss}{Y} %
\begin{document}

\title{Convergence of random series and the rate of convergence of the strong law of large numbers in game-theoretic probability}
\author{Kenshi Miyabe\footnote{Research Institute for Mathematical Sciences,  Kyoto University} \ and 
Akimichi Takemura\footnote{Graduate School of Information Science and Technology,  The University of Tokyo}
}

\date{November 2011}
\maketitle

\begin{abstract}
We give a unified treatment of the convergence of random series and the
rate of convergence of strong law of large numbers in the framework of
game-theoretic probability of Shafer and Vovk \cite{ShaVov01}.  We consider games
with the quadratic hedge as well as more general weaker hedges.  The latter
corresponds to existence of an absolute moment of order smaller than two 
in the measure-theoretic framework. 
We prove some precise relations between the convergence of centered random
series and the convergence of the series of prices of the hedges.
When interpreted in measure-theoretic framework, these results characterize
convergence of a martingale in terms of convergence of the series of
conditional absolute moments.
In order to prove these results we 
derive some fundamental results on deterministic strategies of Reality,
who is a player in a protocol of game-theoretic probability.  
It is of particular interest, since  Reality's strategies
do not have any counterparts in measure-theoretic framework, ant yet
they can be used to prove results which can be interpreted in measure-theoretic
framework.
\end{abstract}

\noindent
{\it Keywords and phrases:} \ 
Kronecker's lemma,
law of the iterated logarithm,
L\'evy's extension of Borel-Cantelli lemma,
Marcinkiewicz-Zygmund strong law,
three-series theorem.

\section{Introduction}
\label{sec:intro}

In standard textbooks on measure-theoretic probability, 
the strong law of large numbers (SLLN) is proved using Kronecker's lemma.
As a precondition for Kronecker's lemma, the convergence of a random series is
usually stated in the form of three-series theorem.  Game-theoretic counterpart 
in Section 4.2 of  Shafer and Vovk (\cite{ShaVov01}) basically follows
the same line of argument. However game-theoretic forms of various conditions
for convergence of random series have not been studied in detail.  Indeed,
game-theoretic counterparts of the standard three-series theorem have 
to be stated more carefully than in the measure-theoretic setting
in several respects, such as the treatment of truncation of random variables and 
the martingale nature of game-theoretic framework.
In particular we need to take into account 
Gilat's counter example (\cite{gilat1971}, see also \cite{brown1971},
\cite{brown1972erratum}) to three-series theorem  for martingales.  

In this paper we give a unified game-theoretic 
treatment of convergence of random series and the
rate of convergence of SLLN.
We consider games with the quadratic hedge.
For an i.i.d.\ case in the terminology of measure-theoretic probability, the
law of the iterated logarithm (LIL) gives the precise rate.
In game-theoretic probability which corresponds to a non-identical case, the rate of convergence may be slower.
We give the precise rate in the game-theoretic framework.

We also consider games with more general weaker hedges.
Marcinkiewicz-Zygmund strong law (\cite{marcinkeiwicz-zygmund},
\cite{gut2005}, \cite{KumTakTak07})
suggests that the rate of
convergence of  random series and SLLN should depend on the
existence of moments. In Section \ref{sec:general-hedge} we will show
that the rate is determined by the inverse function of the hedge function.

In order to derive results on convergence of random series,
we study some topics.  One topic is the set of convergence of a martingale.
If a capital process is required to be non-negative, 
the convergence theorem for a non-negative capital process stated in \cite{ShaVov01}
is sufficient.
However it is useful to consider the set of convergence for an arbitrary capital process which may be negative.

Another topic is deterministic strategies of Reality.
We propose a notion called ``compliance''
concerning deterministic strategies of Reality and prove some results on them.
In \cite{ShaVov01}, Reality's strategy is only 
briefly discussed.  Furthermore  only randomized strategies of Reality are considered
in Section 4.3 of \cite{ShaVov01} and Section 7 of \cite{KumTakTak07}.
Our deterministic strategies of Reality can be understood as  ``derandomizations''
of the randomized strategies in \cite{ShaVov01} and \cite{KumTakTak07}.
It is of interest that deterministic strategies of Reality, which do not
have any counterparts in measure-theoretic probability, can be used to prove
results which can be interpreted in measure-theoretic probability.

The organization of this paper is as follows.
In Section \ref{sec:standard-bounded-game}
we consider sets of convergence in the bounded forecasting game and 
establish preliminary results on the implication of 
the convergence of a random series to 
the convergence of the series of prices for Reality's moves.
We also treat the coin-tossing game as a special case of
the bounded forecasting game.
In Section \ref{sec:bound} we consider bounded 
forecasting game with quadratic hedge and prove various results
on convergence of random series.
In Section \ref{sec:rate} we study the rate of convergence
of SLLN in unbounded forecasting game with quadratic hedge
and 
in Section \ref{sec:general-hedge} we generalize our results 
to games with more general weaker hedges.
We end the paper with some discussion on further topics 
in Section \ref{sec:disc}.

\section{Preliminary results for bounded forecasting game and the coin-tossing game}
\label{sec:standard-bounded-game} 

In this section we consider the bounded forecasting game of Section 3.3
of \cite{ShaVov01}
and the coin-tossing game as a special
case.
This section also serves as a brief
introduction to game-theoretic probability.

\subsection{Definitions and some notions}
\label{subsec:preliminary} 
Consider a perfect information game among three players: Forecaster, Skeptic
and Reality.  Let $C>0$ be given. 
Before the start of the game, Skeptic announces
his initial capital $\cK_0 = D > 0$. 
$\cK_0\equiv 1$ in Section 3.3 of \cite{ShaVov01}, 
but in this section, for our discussion of the bounded forecasting game,
it is more convenient to let Skeptic announce 
his initial capital $\cK_0 = D$.  
Then, at each round $n=1,2,\dots$,  of the game,  these players announce
their moves in the order:  Forecaster, Skeptic
and Reality.  At each round, Forecaster first
announces $m_n$, which is interpreted as the price for Reality's move $x_n$. 
Forecaster has to announce a  {\it coherent} price $m_n$ (Section 1.2 of \cite{ShaVov01}), that is, with the announced price  $m_n$, 
Reality should always be able to prevent Skeptic
from strictly increasing his capital $\cK_n$ for sure.  
Given the price, Skeptic then announces the amount $M_n$  %
he bets. %
Finally Reality announces her move $x_n \in [-C, C]$.
The payoff to Skeptic at the $n$-th round is $M_n (x_n - m_n)$ and his capital
is updated as $\cK_n = \cK_{n-1} + M_n (x_n - m_n)$.

More precisely, the protocol of bounded 
forecasting game %
is written as follows. 

\begin{quote}
{\sc Bounded  Forecasting Game}\\*
\textbf{Parameter}: $C>0$\\*
\textbf{Players}: Forecaster, Skeptic, Reality\\
\textbf{Protocol}:\\
\indent Skeptic announces his initial capital $\cK_0=D>0$.\\
\indent FOR $n=1,2,\ldots$:\\
\indent\indent Forecaster announces $m_n\in[-C,C]$.\\
\indent\indent Skeptic announces $M_n\in\mathbb{R}$.\\
\indent\indent Reality announces $x_n\in[-C,C]$.\\
\indent\indent $\cK_n:=\cK_{n-1}+M_n(x_n-m_n)$.\\
\textbf{Collateral Duties}:
Skeptic must keep $\cK_n$ non-negative.
Reality must keep $\cK_n$ from tending to infinity.
\end{quote}

Note that Forecaster's price  $m_n \in [-C,C]$ is clearly coherent,
because Reality can always choose $x_n = m_n$, so that $\cK_n = \cK_{n-1}$ 
irrespective of $M_n$.

A strategy $\cP=\{\cP_n\}_{n\ge 1}$ of Skeptic specifies $M_n$, $n\ge 1$,  in terms
of past moves of other players $m_k, x_k$, $k=1,\dots,n-1,$ and
the current price $m_n$:
\[
M_n = \cP_n(m_1, x_1,\dots, m_{n-1},x_{n-1},m_n).
\]
We define the capital process $\cK^\cP=\{\cK^{\cP}_n\}_{n\ge 0}$ for a given strategy $\cP$ recursively by $\cK_0^\cP=0$ and 
\[
\cK^{\cP}_n = \cK^{\cP}_{n-1} + \cP_n(m_1, x_1,\dots, m_{n-1},x_{n-1},m_n) (x_n - m_n),
\ \  n\ge 1. 
\]
$\cK_n^{\cP}$ is the cumulative payoff to Skeptic up to round $n$ 
under the strategy $\cP$ (without the initial capital $\cK_0=D$). 
With our definition, $\cK_0^\cP\equiv 0$ is distinguished from the initial capital
$\cK_0=D$ announced by Skeptic.
We call a sequence of real-valued functions $\cS_n(m_1,x_1, \dots, x_n, m_n)$ of
$m_1, x_1, \dots, m_n, x_n$, $n\ge 0$, 
a capital process if $\cS_n=\cK_n^\cP$ for some strategy $\cP$.

An infinite sequence $\xi=(m_1, x_1, m_2, x_2, \dots)$ 
of moves of Forecaster and Reality is called a {\it path}.
Define the sample space 
\[
\Xi=\{ \xi=(m_1, x_1, m_2, x_2, \dots) \mid m_n, x_n \in [-C,C], \forall n\ge 1 \}
\]
as the set of  paths.
We regard $\cK^\cP$ and $\cP$ as functions of $\xi$.
They are denoted by $\cK_n^\cP(\xi)$ and $\cP_n(\xi)$,
which actually depend only on prefixes of $\xi$ of length $2n$ and $2n-1$, respectively.

Any subset $E$ of $\Xi$ is called an event.  We say that
Skeptic can {\it force} $E$ if there exists a strategy $\cP$ of Skeptic, 
such that
\begin{equation}
\label{eq:collateral-duty}
\cK_n^\cP(\xi)\ge -1, \ \ \forall \xi\in \Xi, \forall n\ge 0,
\end{equation}
and 
\begin{equation*}
\xi \not \in E \ \Rightarrow\ 
\limsup_n \cK^{\cP}_n(\xi) = \infty.
\end{equation*}
In this paper we do not make the distinction between
forcing ($\lim_n \cK^{\cP}_n=\infty$) and weak forcing ($\limsup_n \cK^{\cP}_n=\infty$) 
in view of Lemma 3.1 of \cite{ShaVov01}.
A strategy $\cP$ satisfying \eqref{eq:collateral-duty} is called {\it prudent}, i.e.,
if Skeptic observes his collateral using $\cP$ with the initial capital
$\cK_0=1$. 

For two events $E_1, E_2\subset \Xi$, the
event $E_1 \Rightarrow E_2$ stands for $E_1^C \cup E_2$, where
$E_1^C$ is the complement of $E_1$.  $E_1 \Leftrightarrow E_2$
stands for both implications:
\[
(E_1^C \cup E_2)\cap (E_1 \cup E_2^C) = (E_1 \cap E_2) \cup (E_1^C \cap E_2^C)
= (E_1 \Delta E_2)^C,
\]
where $\Delta$ denotes the symmetric difference.  Note that
$E_1 \Leftrightarrow E_2$ and 
$E_1^C \Leftrightarrow E_2^C$ are the same as a subset of $\Xi$.
In  view of Lemma 3.1 of \cite{ShaVov01}, Skeptic can force $E_1 \Leftrightarrow E_2$
if and only if Skeptic can force both $E_1 \Rightarrow E_2$ and $E_2 \Rightarrow E_1$.

\subsection{A set of convergence}
\label{subsec:set-convergence}

Martingale convergence theorems in measure-theoretic probability state that
the limit of a martingale exists and is finite almost surely if the martingale is bounded in $\mathcal{L}^1$.
In Section VII.5 of \cite{Shi95} the set of convergence was studied when the condition is not satisfied.

Game-theoretic probability also has convergence theorems.
If a capital process is required to be non-negative, the convergence theorem holds.
However it is useful to consider a strategy whose capital process may be negative
in order to construct a strategy whose capital process is non-negative as we will do later.
Then we will prove a game-theoretic version of a simple case of results in \cite{Shi95}.
The results are used in a later section.

Let $\cP$ be a strategy of Skeptic.  Denote the
minimum possible gain $L_n=L_n(m_1, x_1, \dots, m_{n-1}, x_{n-1}, m_n)$
to Skeptic at the round $n$ (after he knows Forecaster's move $m_n$)
under $\cP$ by
\begin{align*}
L_n  &= \min_{x\in [-C,C]} \cP_n(m_1, x_1,\dots,m_{n-1}, x_{n-1}, m_n)(x-m_n)\\
     &= \begin{cases} \cP_n\times (-C-m_n) & \text{if } 
 \cP_n\ge 0 \\
 \cP_n \times (C-m_n) & \text{otherwise}.
\end{cases}
\end{align*}
For $D>0$, we define 
the stopping time $\tau_\cP^D=\tau_\cP^D(\xi)$ as the first time $D+ \cK_n^\cP$ may be negative:
\begin{equation*}
\tau_\cP^D=\min\{n \ge 1: \cK_1^\cP\ge -D, \dots, \cK_{n-1}^\cP \ge -D,\  \cK_{n-1}^\cP + L_n  < -D\}.
\end{equation*}
As usual $\tau_\cP^D=+\infty$ if the set on the right-hand side is empty.
The truncation $\cP^D$ of $\cP$ at the loss $-D$ is defined as
\[
\cP^D_n(m_1, x_1,\dots, m_n)
= \begin{cases} \cP_n(m_1, x_1,\dots, m_n) & \text{if } n < \tau^D_\cP\\
  0 & \text{otherwise}.
\end{cases}
\]
Note that starting with the initial capital of $D>0$, 
Skeptic observes his collateral duty by employing $\cP^D$, i.e.,
$D+\cK_n^{\cP^D}$ is always nonnegative. 

We now prove the following proposition.

\begin{proposition}
\label{prop:version-of-martingale-convergence}
Let $\cP$ be any strategy in the bounded forecasting game. Let
$B^\cP$ denote the event 
\begin{equation}
\label{eq:bounded-bet}
B^\cP = \{ \cP_n \text{ is bounded }\}
=\{\xi \mid \sup_n | \cP_n(\xi)| < \infty\}.
\end{equation}
Skeptic can force
\begin{equation}
\label{eq:kp-convergence}
B^\cP \ \Rightarrow\ \big(\cK^\cP_n \text{ converges in $\mathbb R$  \ or }\ 
 (\limsup_n \cK^\cP_n =+\infty \ \text{ and } \ 
\liminf_n \cK^\cP_n = -\infty) \big).
\end{equation}
\end{proposition}
By ``convergence in $\mathbb R$''
we mean that  $\lim_n \cK_n^\cP$ 
exists and is finite.  In later statements we will omit ``in $\mathbb R$'' for
simplicity.  Proposition \ref{prop:version-of-martingale-convergence}
means that given any strategy  $\cP$, there exists another prudent strategy
$\cQ$ of Skeptic, such that $\limsup_n \cK_n^\cQ=\infty$ 
if \eqref{eq:kp-convergence} is violated.

\begin{proof}
Note that the convergence or divergence of $\cK_n^\cP$
is classified into the following five exclusive
cases: 
\begin{itemize}
\setlength{\itemsep}{0pt}
\item[(i)] $-\infty < \liminf_n \cK_n^\cP = \limsup_n \cK_n^\cP < \infty$ \quad (convergence),
\item[(ii)] $ -\infty < \liminf_n \cK_n^\cP <  \limsup_n \cK_n^\cP < \infty$ \quad (bounded oscillation),
\item[(iii)] $-\infty <  \liminf_n \cK_n^\cP \le  \limsup_n \cK_n^\cP = \infty$, 
\item[(iv)]  $-\infty =  \liminf_n \cK_n^\cP \le  \limsup_n \cK_n^\cP < \infty$, 
\item[(v)] $-\infty = \liminf_n \cK_n^\cP, \ \limsup_n \cK_n^\cP = \infty$ \  (two-sided unbounded oscillation).
\end{itemize}
According to this classification, the sample space $\Xi$ is partitioned
into five subsets $E_1^\cP, \dots, E_5^\cP$.
By Lemma 3.2 of \cite{ShaVov01},
it suffices to construct a prudent strategy $\cQ$  of Skeptic
for each of the cases $E_i^\cP$, $i=2,3,4$, such that
$\xi\in B^\cP \cap E_i^\cP$ implies  $\limsup_n \cK_n^\cQ=\infty$, $i=2,3,4$.

We consider the case (iii) in detail.  
As noted above, for each $D>0$,  $D + \cK_n^{\cP^D}$ is always
nonnegative.  Consider dividing the initial capital of 1 into 
countably infinite accounts with initial capitals $1/2 + 1/4 + \dots = 1$.
For the $D$-th account with the initial capital of $1/2^D$, we apply the strategy 
$\cP^D/(D 2^D)$.  The resulting combined strategy 
$\cQ$ is written as
\[
\cQ= \sum_{D=1}^\infty \frac{1}{D2^D} \cP^D.
\]
Then the capital process of $\cQ$ is written as
\[
 1 + \cK_n^{\cQ} 
= \sum_{D=1}^\infty \frac{1}{2^D} 
+ \sum_{D=1}^\infty \frac{1}{D2^D} \cK_n^{\cP^D}
 =\sum_{D=1}^\infty \frac{1}{D2^D} (D + \cK_n^{\cP^D})
\]
and hence $\cQ$ is prudent. Now for each $\xi\in B^\cP\cap E_3$,
there exist positive constants 
$ D_1= D_1(\xi), D_2 = D_2(\xi)$, such that
\[
 \cK_n^\cP(\xi)> -D_1(\xi),  \quad  |\cP_n(\xi)| < D_2(\xi),  \quad \forall n\ge 1.
\]
Then since $|m_n|, |x_n| \le C$
\[
|L_n(\xi)| \le 2 C D_2(\xi), \ \forall n\ge 1.
\]
Consider $D > D_1(\xi) + 2 C D_2(\xi)$.  Then for all $n\ge 1$,
we have 
\[
\cK_k^\cP(\xi)\ge -D_1(\xi) > -D,\ \  k=1,\dots,n-1,
\]
and
\[
\cK_{n-1}^{\cP}(\xi) + L_n(\xi) \ge -D_1(\xi)-2C D_2(\xi) > -D.
\]
Hence for this $D$ we have $\tau_\cP^D(\xi)=\infty$ and
$
\cK_n^\cP=\cK_n^{\cP^D},  \ \forall n\ge 1.
$
Therefore 
$\xi\in B^{\cP} \cap E_3\ \Rightarrow \ \limsup_n \cK_n^{\cQ}(\xi)=\infty$.  
This proves the case of (iii).

The case (iv) is proved by the symmetry of the bounded
forecasting protocol, i.e.\ by considering $-\cP$ instead of $\cP$.

Finally  (ii) can be proved by the standard argument
involving Doob's upcrossing lemma (see Lemma 4.5 of \cite{ShaVov01}).
Note that Lemma 4.5 of \cite{ShaVov01} is for the case of
prudent $\cP$.  In our case, $\cP$ is not necessarily
prudent.  However again combining truncations $\cP^D$, $D=1,2,\dots$,
with the argument of upcrossing lemma, 
we can construct a prudent $\cQ$ such that $\limsup_n \cK_n^{\cQ}(\xi)=\infty$
for each $\xi\in B^\cP\cap E_2$.  This proves the proposition.
\end{proof}

We call a strategy $\cP$ {\it uniformly cautious} if 
\[
\sup_{\xi\in \Xi, n\ge 1} |\cP_n(\xi)| < \infty.
\]
For uniformly cautious $\cP$,  $B^\cP=\Xi$.
Therefore if $\cP$ is uniformly cautious, then Skeptic
can force the right-hand side of  \eqref{eq:kp-convergence}.
The reason for considering $B^\cP$ in \eqref{eq:bounded-bet}
and uniformly cautious strategies  is that they eliminate
doubling type strategies.

Using Proposition \ref{prop:version-of-martingale-convergence}
we can prove the following result, which is 
a generalization of results in Section 2.2.2 of \cite{horikoshi-takemura-2008}.

\begin{proposition}\label{prop:bound-ht}
In the bounded forecasting game Skeptic can force
\begin{align}
\sum_n x_n \text{ converges}\ \Rightarrow \ \   
&  \sum_n m_n\text{ converges  \ \ or}\nonumber\\
& (\limsup_n \sum_{k=1}^n  m_k = +\infty \ \text{ and } \ 
\liminf_n \sum_{k=1}^n  m_k = -\infty)
\label{eq:xnconvergence}
\end{align}
and
\begin{align}
\sum_n m_n \text{ converges}\ \Rightarrow \ \   
&  \sum_n x_n\text{ converges \ \ or}\nonumber\\
& (\limsup_n \sum_{k=1}^n  x_k = +\infty \ \text{ and } \ 
\liminf_n \sum_{k=1}^n  x_k = -\infty)
\label{eq:mnconvergence}
\end{align}
\end{proposition}

\begin{proof}
Let $\ss_n = \sum_{k=1}^n (x_k - m_k)$.
Consider the uniformly cautious strategy $\cP$ such that $M_n\equiv 1$.  Then 
$\cK_n^{\cP}=\ss_n$.
By Proposition \ref{prop:version-of-martingale-convergence}, 
Skeptic can force 
\[
\ss_n \text{ converges  or }\ 
 (\limsup_n \ss_n = +\infty \ \text{ and } \ 
\liminf_n \ss_n  = -\infty).
\]
We separate the sums in $\ss_n$ as $\ss_n = \sum_{k=1}^n x_k - \sum_{k=1}^n m_k$ and
restrict relevant events to the  particular 
event $E_0=\{\sum_n x_n\ \text{converges}\}$. 
Then clearly Skeptic can force \eqref{eq:xnconvergence}.
\eqref{eq:mnconvergence} is proved similarly, by switching 
the roles of $x_n$ and $m_n$.
\end{proof}

\subsection{Some applications}
\label{subsec:application}

Consider the multi-dimensional bounded forecasting game
(cf.\ \cite{KumTakTak11}) defined as follows.
Reality's move space is a compact set $\cal X$ of ${\mathbb R}^d$.
Forecaster's move space  is
the convex hull ${\rm co}({\cal X})$ of $\cal X$
and Skeptic's move space is 
${\mathbb R}^d$.  
Denote the moves by Forecaster, Skeptic and Reality by
$\mu_n$, ${\mathfrak m}_n$  and $\chi_n$, respectively.
The payoff to Skeptic is ${\mathfrak m}_n \cdot (\chi_n - \mu_n)$, where
``$\cdot$'' denotes the standard inner product in ${\mathbb R}^d$.
The protocol of the multi-dimensional bounded forecasting game
is written as follows.

{\rm
\begin{quote}
  {\sc Multi-dimensional Bounded Forecasting Game} \\
\textbf{Parameter}: a compact region ${\cal X}\subset {\mathbb R}^d$\\
\textbf{Players}: Forecaster, Skeptic, Reality\\
\textbf{Protocol}:\\
\indent Skeptic announces his initial capital $\cK_0=D>0$.\\
\indent FOR $n=1,2,\ldots$:\\
\indent\indent Forecaster announces $\mu_n\in {\rm co}{\cal X}$.\\
\indent\indent Skeptic announces ${\mathfrak m}_n\in\mathbb{R}^d$.\\
\indent\indent Reality announces $\chi_n\in {\cal X}$.\\
\indent\indent $\cK_n:=\cK_{n-1}+{\mathfrak m}_n\cdot (\chi_n-\mu_n)$.\\
\textbf{Collateral Duties}:
Skeptic must keep $\cK_n$ non-negative.
Reality must keep $\cK_n$ from tending to infinity.
\end{quote}
}
\vspace*{-3pt}
It is easily seen that 
Proposition \ref{prop:version-of-martingale-convergence}
holds in the multi-dimensional bounded forecasting game.

We now consider the  coin-tossing game and give a 
game-theoretic probability version of L\'evy's extension of Borel-Cantelli lemma.
The protocol of the coin-tossing game is written as follows.
\begin{quote}
{\sc Coin-tossing game}\\
\textbf{Protocol}:\\
\indent $\cK_0=1$.\\
\indent FOR $n=1,2,\ldots$:\\
\indent\indent Forecaster announces $p_n\in[0,1]$\\
\indent\indent Skeptic announces $M_n\in\mathbb{R}$.\\
\indent\indent Reality announces $x_n\in \{0,1\}$.\\
\indent\indent $\cK_n:=\cK_{n-1}+M_n(x_n-p_n)$.\\
\end{quote}

We have the following. 
\begin{example} (A game-theoretic version of L\'evy's extension of Borel-Cantelli lemma)\quad 
\label{th:coin}
In the coin-tossing game
Skeptic can force
\begin{equation}
\label{eq:coin-equivalence}
\sum_n p_n<\infty\iff \sum_n x_n<\infty.
\end{equation}
\end{example}

\begin{remark}
The statement may be easier to understand if we rewrite it as follows:
Skeptic can force
\[\sum_n p_n<\infty\Rightarrow \sum_n x_n<\infty\]
and
\[\sum_n p_n=\infty\Rightarrow \sum_n x_n=\infty.\]
\end{remark}

Coin-tossing game is a special case of the bounded forecasting game (with $C=D=1$), 
in such a way that the move space of Reality  is restricted to $\{0,1\}$
and  the move space of Forecaster is restricted to $[0,1]$ by coherence.
Therefore, if Skeptic can force an event $E$ in the bounded
forecasting game, then Skeptic can force $E$ in the coin-tossing game.
In the coin-tossing game, $\sum_n x_n$ and $\sum_n p_n$ are
non-negative series and they either converge to finite values or
diverge to $+\infty$.  Therefore the case of two-sided unbounded 
oscillation on the right-hand side of 
\eqref{eq:xnconvergence} and
\eqref{eq:mnconvergence} is impossible, which implies that Skeptic can force
\eqref{eq:coin-equivalence}.

L\'evy's extension of Borel-Cantelli lemma in measure-theoretic
probability is usually stated as follows (cf.\ Theorem 12.15 of \cite{Wil91}).
\begin{proposition}\label{pro:levy}
Let $X_n$, $n=1,2,\dots,$ be a sequence of $0$-$1$ random
variables adapted to filtration $\{{\cal F}_n\}$.  Let $p_n=E(X_n\mid {\cal F}_{n-1})$.
Then almost surely 
\begin{enumerate}
\item $\sum_n p_n < \infty \ \Rightarrow\ \sum_n x_n < \infty$,
\item $\sum_n p_n = \infty \ \Rightarrow \ \lim_N (\sum_{n=1}^N x_n / \sum_{n=1}^N p_n) = 1$.
\end{enumerate}
\end{proposition}
If we weaken (ii) to $\sum_n p_n =\infty \Rightarrow \sum_n x_n = \infty$, then
the measure-theoretic extension looks similar to 
Theorem \ref{th:coin}.  However there are some basic differences.
In our setting, $p_n$'s are only ``prequentially'' 
(eg.\ \cite{vovk-shen-2010}) announced by Forecaster and
there is no need to specify the full probability measure on
$x_n$, $n=1,2,\dots$.  Also, in measure-theoretic framework 
the null set, 
where these implications do not hold, may depend on the underlying probability measure.
On the other hand, in the game-theoretic setting, 
we have constructed an explicit strategy $\cQ$ forcing 
\eqref{eq:coin-equivalence} and the behavior of
its capital process $\cK^\cQ$  on the symmetric difference of two sets in \eqref{eq:coin-equivalence} is explicitly understood.
Furthermore in Proposition \ref{prop:levy-ii} of Section \ref{sec:rate} 
we will strengthen the rate of convergence in (ii).

\section{Bounded forecasting game with quadratic hedge}
\label{sec:bound}

The standard measure-theoretic three-series theorem  involves truncation of 
random variables and their means and variances.  In considering game-theoretic
counterpart of the standard setup, we here  consider the simple
case that the truncation is given before the game,
i.e.\ we 
consider a variant of bounded forecasting game.
In addition we assume that   
the quadratic hedge is available to Skeptic.  
From now on for simplicity we assume $\cK_0=1$.
The protocol for this section is written as follows.

\begin{quote}
{\sc Bounded Forecasting with Quadratic Hedge (BFQH)}\\
\textbf{Parameter}: $C>0$\\
\textbf{Players}: Forecaster, Skeptic, Reality\\
\textbf{Protocol}:\\
\indent $\cK_0=1$.\\
\indent FOR $n=1,2,\ldots$:\\
\indent\indent Forecaster announces $m_n\in[-C,C]$ and $v_n\in[0,C^2-m_n^2]$.\\
\indent\indent Skeptic announces $M_n\in\mathbb{R}$ and $V_n\in\mathbb{R}$.\\
\indent\indent Reality announces $x_n\in[-C,C]$.\\
\indent\indent $\cK_n:=\cK_{n-1}+M_n(x_n-m_n)+V_n((x_n-m_n)^2-v_n)$.\\
\textbf{Collateral Duties}:
Skeptic must keep $\cK_n$ non-negative.
Reality must keep $\cK_n$ from tending to infinity.
\end{quote}

In this protocol the following theorem holds.

\begin{theorem}\label{th:bound}
In BFQH Skeptic can force
\[\sum_n v_n<\infty\iff\sum_n(x_n-m_n)\textrm{ converges.}\]
\end{theorem}

Before giving a proof of this theorem, we discuss some features
of the theorem and the protocol BFQH.

The implication ``$\sum_n v_n < \infty \ \Rightarrow\ \sum_n (x_n - m_n) \text{ converges}$'' holds even for the case of $C=+\infty$, but 
the the converse implication does not hold for $C=+\infty$,
as shown in Lemma \ref{lem:vconv} below.  Therefore
the main point of Theorem \ref{th:bound} is that Skeptic can force
the equivalence of both sides in the case $0 <  C < +\infty$.

In BFQH the ranges of $V_n$ and $v_n$ are different from %
both the bounded forecasting game in Section 3.1 of \cite{ShaVov01}
and the unbounded forecasting game in Section 4.1 of \cite{ShaVov01}.
First, we allow $V_n$ to be negative.
In the usual unbounded protocol, 
if Skeptic announces negative $V_n$ then it violates his collateral duty because $\cK_n\to-\infty$ as $x_n\to\infty$.
However this does not happen in the above protocol because $|x_n|\le C$ is bounded.
Second, we restrict $v_n\in[0,C^2-m_n^2]$ for the game to be coherent.
For example, consider the case $m_n=0$, $v_n=C^2+\epsilon$ $(\epsilon>0)$, $M_n=0$ and $V_n=-1$,
Then $\cK_n=\cK_{n-1}-x_n^2+C^2+\epsilon\ge \cK_{n-1}+\epsilon$ for all $x_n\in[-C,C]$.
More precisely we should restrict $m_n$ and $v_n$ such that 
\begin{align}\label{eq:restmv}
m_n=\int_{[-C,C]}zdp_n\textrm{ and }v_n=\int_{[-C,C]}(z-m_n)^2dp_n
\end{align}
for some probability measure $p_n$.
For $v_n\in[0,C^2-m_n^2]$ it is easily checked that
(\ref{eq:restmv}) holds
for the two-point measure
\[p_n(\{a_n\})=\frac{1-m_n/C}{2}\mbox{ and }p_n(\{b_n\})=\frac{1+m_n/C}{2},\]
where $a_n=m_n+(-C-m_n)\sqrt{\frac{v_n}{C^2-m_n^2}}$ and 
$b_n=m_n+(C-m_n)\sqrt{\frac{v_n}{C^2-m_n^2}}$.

BFQH can be regarded as a variant 
of the two-dimensional bounded forecasting game 
in Section \ref{subsec:application}
with an  additional restriction at each round to Reality.
In the multi-dimensional bounded forecasting game let
${\cal X}=[-2C,C]\times [0,4C^2]$,
$\mu_n =(m_n, v_n)$, $\fm_n=(M_n, V_n)$,
$\chi_n=(x_n , (x_n-m_n)^2)$. Also at each round of the game
put an additional restriction to Reality's move space depending on Forecaster's move
as $\chi_n \in {\cal X}_n = \{(x,(x-m_n)^2) | -C \le x\le C\}$.
Since the restriction is advantageous to Skeptic, Proposition
\ref{prop:version-of-martingale-convergence} holds for BFQH.

Now we give a proof of Theorem  \ref{th:bound}. %
We use the notation 
\[
\ss_n = \sum_{k=1}^n (x_k - m_k), \ \ss_0=0,  \ \ A_n = \sum_{k=1}^n v_k,  \ A_\infty = \sum_{k=1}^\infty v_k.
\]

\begin{proof}[Proof of Theorem \ref{th:bound}]
\noindent($\Rightarrow$)
Consider a capital process 
\[
{\mathcal T}_n = \ss_n^2 - A_n = 2 \sum_{k=1}^n \ss_{k-1} (x_k-m_k)  + \sum_{k=1}^n ((x_k - m_k)^2 - v_k),
\]
for Skeptic's strategy $M_n = 2 \ss_{n-1}$ and $V_n = 1$. $A_n$ is the
compensator for $\ss_n^2$.   Then by Lemma 4.6 and Lemma 4.7 of \cite{ShaVov01},
Skeptic can force that $\ss_n$ converges.

\medskip
\noindent
($\Leftarrow$)
Although the argument for this implication is essentially the 
same as the first part, we can not directly apply 
Lemma 4.6 of \cite{ShaVov01}.  We prove this implication
by Proposition \ref{prop:version-of-martingale-convergence}.
Consider Skeptic's  strategy 
\begin{equation*}
\cP_0: 
M_n= -2 \ss_{n-1}, \ V_n =-1\le0,
\end{equation*}
which is the negative of the above strategy. 
The capital process of $\cP_0$ is given as
$\cK_n^{\cP_0}=A_n - \ss_n^2$.  %
Note that $B^{\cP_0}$ of Proposition \ref{prop:version-of-martingale-convergence}
for this strategy is the set of paths such that $\{\ss_n\}$ is bounded.
Therefore by a multi-dimensional version of Proposition \ref{prop:version-of-martingale-convergence},
Skeptic can force 
\begin{align}
\{Y_n\} \textrm{ is bounded} & \ \Rightarrow \ 
\big ( A_n - \ss_n^2 \text{ converges}  \text{\  or\ }\nonumber\\
& \qquad (\limsup_n (A_n - \ss_n^2) = +\infty \text{ and }
\liminf_n (A_n - \ss_n^2) = -\infty) \big).
\label{eq:an-reverse-implication}
\end{align}
By assumption $\ss_n$ converges. 
Then $\{\ss_n\}$ is bounded. Also $\ss_n^2$ converges.
Furthermore since  $A_n$ is non-decreasing, the second case of 
the right-hand side of 
(\ref{eq:an-reverse-implication}) is impossible.
Therefore Skeptic can force the event: $\ss_n$ converges $\Rightarrow A_\infty < \infty$.
\end{proof}

In Theorem \ref{th:bound} we considered the convergence of the series
$\sum_n(x_n-m_n)$.
In standard measure-theoretic probability theory, in the bounded case,
the series is split as $\sum_n (x_n-m_n) = \sum_n x_n - \sum_n m_n$.
Then the convergence of $\sum_n x_n$ is discussed in terms of
convergence of $\sum_n m_n$ and $\sum_n v_n$.
Two-series theorem (see e.g.\ \cite{Shi95}) says that,
under the assumption of independence,
$\sum_n x_n$ converges if and only if
$\sum_n v_n$ converges and $\sum_n m_n$ converges.  
However, without the assumption of independence,
the case where $\sum_n m_n$ does not converge and $\sum_n v_n=\infty$,
is very subtle. 
Indeed, Gilat \cite{gilat1971} 
showed a counter example which implies that there exists no condition on the conditional expectations
and conditional variances of a series of uniformly bounded random variables
which is necessary and sufficient for the sum to converge almost surely.
We give a game-theoretic version of Gilat's result.

\begin{theorem}
In BFQH there is no set $D$ of infinite sequences of Forecaster's moves such that Skeptic can force
\begin{align}\label{th:convset}
(m_1,v_1,m_2,v_2,\ldots)\in D\iff \sum_n x_n 
\mbox{ converges. }
\end{align}
\end{theorem}

First we explain the intuition behind the proof.
The following is Gilat's example.
Let $\{Z_n\}$ be the sequence of independent random variables with
\[P\{Z_n=n^{-1/2}\}=\frac{1}{2}=P\{Z_n=-n^{1/2}\}\mbox{ for }n\ge1,\]
and $Z_0=0$.
Then $Z_n\to0$ and $\sum_n Z_n$ does not converge almost surely.
Let $X_n=Z_n-Z_{n-1}$ and $Y_n=-Z_n-Z_{n-1}$ for $n\ge1$.
Then the conditional means and the conditional variances of $X_n$
are the same as those of $Y_n$.
Furthermore $\sum_{k=1}^n X_k=Z_n\to0$ and $\sum_{k=1}^n  Y_k= -2 \sum_{k=1}^{n-1} Z_k - Z_n$ does not converge almost surely.

We want to show that a strategy forcing \eqref{th:convset} yields 
a contradiction.
Intuitively, if Skeptic can force an event  $E$, then a ``typical'' or ``random'' sequence should satisfy $E$.
Pick up a typical realization  $\{\eta_n\}$ of $\{Z_n\}$.
Then $\eta_n\to0$ and $\sum_n \eta_n$ does not converge.
The sequence $\{x_n\}=\{\eta_n-\eta_{n-1}\}$ is also typical  for $\{X_n\}$.
Then $\sum_{k=1}^n x_k=\eta_n\to0$ and $(m_1,v_1,m_2,v_2,\ldots)\in D$.
Furthermore $\{y_n\}=\{-\eta_n-\eta_{n-1}\}$ is typical for $\{Y_n\}$.
Then $\sum_n y_n$ does not converge and $(m_1,v_1,m_2,v_2,\ldots)\not\in D$.
This is a contradiction.
Since Forecaster knows the previous moves of Reality,
the first randomness does not imply the second and third randomness
in game-theoretic probability.
Hence we construct three strategies and pick up a random sequence
for the strategies.

\begin{proof}
Assume that there exists a set $D$ and a (prudent) strategy $\cQ$
by which Skeptic can force that
\begin{align}\label{eq:Qforce}
(m_1,v_1,m_2,v_2,\ldots)\in D\iff S_n\mbox{ converges. }
\end{align}
We prove that this assumption yields a contradiction.

Also fix a (prudent) strategy $\cP$ of Skeptic that forces
\begin{align}\label{eq:Pforce}
\sum_n v_n<\infty\iff\sum_n(x_n-m_n)\mbox{ converges.}
\end{align}
Such a strategy exists by Theorem \ref{th:bound}.

Let $a\in\mathbb{N}$ be a positive constant large enough to satisfy
$0<(a-1)^{-1/2}+a^{-1/2}\le C$
and define
\[\omega_n=(n+a)^{-1}, \quad n\ge 1.
\]
Then $\omega_{n-1}^{1/2}+\omega_n^{1/2}<(a-1)^{-1/2}+a^{-1/2}\le C$.
Hence
\[\omega_{n-1}^{1/2}+\omega_n^{1/2}<C\mbox{ and }\omega_{n-1}+\omega_n<C^2.\]

We construct two strategies $\cQ_1$ and $\cQ_2$ for Skeptic from $\cQ$.
At round $n$, Skeptic by $\cQ$ announces $M_n^{\cQ}$ and $V_n^{\cQ}$.
These values depend on the  moves of Forecaster and Reality until
the previous round and the moves of Forecaster at the current round.
So we write
\begin{align*}
M_n^{\cQ}=M_n^{\cQ}(m_1,v_1,x_1,\ldots,m_{n-1},v_{n-1},x_{n-1},m_n,v_n),\\
V_n^{\cQ}=V_n^{\cQ}(m_1,v_1,x_1,\ldots,m_{n-1},v_{n-1},x_{n-1},m_n,v_n).
\end{align*}
Define  $\cQ_1$ and $\cQ_2$ by 
\begin{align*}
M_n^{\cQ_1}=&M_n^{\cQ}(-\eta_0,\omega_1,\eta_1-\eta_0,-\eta_1,\omega_2,\eta_2-\eta_1,\ldots,-\eta_{n-1},\omega_n),\\
V_n^{\cQ_1}=&V_n^{\cQ}(-\eta_0,\omega_1,\eta_1-\eta_0,-\eta_1,\omega_2,\eta_2-\eta_1,\ldots,-\eta_{n-1},\omega_n),\\
M_n^{\cQ_2}=&-M_n^{\cQ}(-\eta_0,\omega_1,-\eta_1-\eta_0,-\eta_1,\omega_2,-\eta_2-\eta_2,\ldots,-\eta_{n-1},\omega_n),\\ 
V_n^{\cQ_2}=&V_n^{\cQ}(-\eta_0,\omega_1,-\eta_1-\eta_0,-\eta_1,\omega_2,-\eta_2-\eta_2,\ldots,-\eta_{n-1},\omega_n),
\end{align*}
where $\eta_0=0$ and $|\eta_n|=\omega_n^{1/2}$. The sign of $\eta_n$
will be defined below in \eqref{eq:eta-def}.
Since the sign of $\eta_n$ is not defined yet, $\cQ_1$ and $\cQ_2$ are not defined yet.
However, we have $-\eta_{n-1}\in[-C,C]$, $\omega_n\in[0,C^2-\eta_{n-1}^2]=[0,C^2-\omega_{n-1}]$ and
$|\pm\eta_n-\eta_{n-1}|\le\eta_n+\eta_{n-1}=\omega_n^{1/2}+\omega_{n-1}^{1/2}\le C$.
Hence, 
if the sign of $\eta_k$ is defined for each $k<n$,  then
$M_n^{\cQ_1}$, $V_n^{\cQ_1}$, $M_n^{\cQ_2}$ and $V_n^{\cQ_2}$ are defined,  
because $M_n^{\cQ}$ and $V_n^{\cQ}$ are defined for
$m_n\in[-C,C]$, $v_n\in[0,C^2-m_n^2]$ and $x_n\in[-C,C]$.

The strategies of $\cQ_1$ and $\cQ_2$ are explained as follows.
$\cQ_1$ moves as $\cQ$ moves when Forecaster's moves are specified as
\begin{equation*}
m_n=-\eta_{n-1},\ v_n=\omega_n
\end{equation*}
and Reality announces $\eta_n-\eta_{n-1}$.
Similarly $\cQ_2$ moves as $\cQ$ when $m_n=-\eta_{n-1}$, $v_n=\omega_n$
and Reality announces $-\eta_n-\eta_{n-1}$
with the minus sign in the definition of $M_n^{\cQ_2}$.

Define a Skeptic's strategy $\cT$ by $\cT=(\cP+\cQ_1+\cQ_2)/3$ and  
define (the sign of) $\eta_n$, $n\ge 1$, by
\begin{equation}
\label{eq:eta-def}
\eta_n=\begin{cases}
\omega_n^{1/2}
&\mbox{ if }M_n^{\cT}\le0\\
-\omega_n^{1/2}
&\mbox{ otherwise.}\end{cases}
\end{equation}
We claim that (for a given $\cP$)
$\cT$ and $\eta_n$ are defined by mutual induction.
Note that $\eta_0$ is already defined as 0.
At round $n$, $\eta_k$ has been defined for each $k<n$.
Then $M_n^{\cQ_1}$ and $M_n^{\cQ_2}$ are defined.
Hence $M_n^{\cT}$ is defined.
It follows that $\eta_n$ is defined.
By induction $\cT$ and $\eta_n$ are defined.

The capital process $\cK_n$ is completely determined by strategies for Forecaster, Skeptic and Reality.
By $\cK_n^{\cF,\cS,\cR}$ we denote the capital at round $n$ when Forecaster uses the strategy $\cF$,
Skeptic $\cS$ and Reality $\cR$.
In this proof we use the following strategies for Forecaster:
\begin{align*}
&\cF_1: m_n=0,\ v_n=\omega_n,\\
&\cF_2: m_n=-\eta_{n-1},\ v_n=\omega_n,
\end{align*}
and for Reality:
\begin{align*}
&\cR_1: x_n=\eta_n,\\
&\cR_2: x_n=\eta_n-\eta_{n-1},\\
&\cR_3: x_n=-\eta_n-\eta_{n-1}.
\end{align*}
We claim $m_n\in[-C,C]$, $v_n\in[0,C^2-m_n^2]$ and $x_n\in[-C,C]$ for each case.
For $\cF_1$, $v_n=\omega_n<C^2$ and $v_n\in[0,C^2]=[0,C^2-m_n^2]$.
For $\cF_2$, $m_n^2=\eta_{n-1}^2=\omega_{n-1}<C^2$ and $m_n\in[-C,C]$,
$v_n=\omega_n<C^2-\omega_{n-1}=C^2-m_n^2$ and $v_n\in[0,C^2-m_n^2]$.
For $\cR_1$, $x_n^2=\eta_n^2=\omega_n<C^2$ and $x_n\in[-C,C]$.
For $\cR_2$ and $\cR_3$, $x_n^2\le(\eta_n+\eta_{n-1})^2\le((a-1)^{-1/2}+a^{-1/2})^2\le C^2$ and $x_n\in[-C,C]$.

\smallskip
\noindent\textbf{Case I.}\quad
Consider the case that Forecaster uses the strategy $\cF_1$, Skeptic uses $\cT$,
and Reality uses $\cR_1$.

Consider the capital $\cK_n^{\cF_1,\cT,\cR_1}$.
By the definition of $\eta_n$ in \eqref{eq:eta-def},
$M_n^{\cT} \eta_n\le 0$ and $\eta_n^2 - \omega_n = 0$ for all $n$.
Hence 
\[\cK_n^{\cF_1,\cT,\cR_1}=\cK_{n-1}^{\cF_1,\cT,\cR_1}+M_n^{\cT}\eta_n+V_n^{\cT}(\eta^2-\omega_n)
\le\cK_{n-1}^{\cF_1,\cT,\cR_1}\]
for all $n$ and $\cK_n^{\cF_1,\cT,\cR_1}$ is monotone non-increasing in this case.
Then
\begin{equation}
\label{eq:KT-bounded}
\cK_n^{\cF_1,\cP,\cR_1}+\cK_n^{\cF_1,\cQ_1,\cR_1}+\cK_n^{\cF_1,\cQ_2,\cR_1} \le 3.
\end{equation}

We will prove
\begin{align}\label{eq:capital-equal}
\cK_n^{\cF_1,\cQ_1,\cR_1}=\cK_n^{\cF_2,\cQ,\cR_2}\mbox{ and }\cK_n^{\cF_1,\cQ_2,\cR_1}=\cK_n^{\cF_2,\cQ,\cR_3}
\end{align}
for all $n$.
Then by prudence of $\cP$ and $\cQ$,
the three capitals $\cK_n^{\cF_1,\cP,\cR_1}$, $\cK_n^{\cF_1,\cQ_1,\cR_1}$ and $\cK_n^{\cF_1,\cQ_2,\cR_1}$ are all non-negative.
Hence they are bounded by \eqref{eq:KT-bounded}.

When Forecaster uses $\cF_1$, $\sum_n m_n=0$ and $\sum_n v_n=\infty$.
By the property (\ref{eq:Pforce}) of $\cP$ and the boundedness of $\cK_n^{\cF_1,\cP,\cR_1}$,
$\sum_n(\eta_n-m_n)=\sum_n \eta_n$ does not converge.
Note that $\eta_n$ converges to 0.

\smallskip
\noindent\textbf{Case II.}\quad
Consider the case that Forecaster uses $\cF_2$, Skeptic $\cQ$ and Reality $\cR_2$.
In this case
\[S_n=\sum_{k=1}^n x_n=\sum_{k=1}^n(\eta_k-\eta_{k-1})=\eta_n,\]
hence it converges.
By the property (\ref{eq:Qforce}) of $\cQ$ and the boundedness of $\cK_n^{\cF_2,\cQ,\cR_2}$,
we have $(-\eta_0,\omega_1,-\eta_1,\omega_2,\ldots)\in D$.

\smallskip
\noindent\textbf{Case III.}\quad
Consider the case that Forecaster uses $\cF_2$, Skeptic $\cQ$ and Reality $\cR_3$.
In this case
\[S_n=\sum_{k=1}^n x_n=\sum_{k=1}^n(-\eta_k-\eta_{k-1})=-2(\sum_{k=1}^n \eta_k)+\eta_n\]
hence it does not converge.
Similarly, by the property (\ref{eq:Qforce}) of $\cQ$ and the boundedness of $\cK_n^{\cF_2,\cQ,\cR_3}$,
we have $(-\eta_0,\omega_1,-\eta_1,\omega_2,\ldots)\not\in D$.
This is a contradiction.

\smallskip
It remains to show the equations (\ref{eq:capital-equal}), that is,
\[\cK_n^{\cF_1,\cQ_1,\cR_1}=\cK_n^{\cF_2,\cQ,\cR_2}\mbox{ and }\cK_n^{\cF_1,\cQ_2,\cR_1}=\cK_n^{\cF_2,\cQ,\cR_3}\]
for all $n$.
We prove this by induction.
For $n=0$, we have $\cK_n^{\cF_1,\cQ_1,\cR_1}=\cK_n^{\cF_2,\cQ,\cR_2}=\cK_n^{\cF_1,\cQ_2,\cR_1}=\cK_n^{\cF_2,\cQ,\cR_3}=1$.
Suppose $\cK_{n-1}^{\cF_1,\cQ_1,\cR_1}=\cK_{n-1}^{\cF_2,\cQ,\cR_2}$ and $\cK_{n-1}^{\cF_1,\cQ_2,\cR_1}=\cK_{n-1}^{\cF_2,\cQ,\cR_3}$.
Then
\begin{align*}
\cK_n^{\cF_1,\cQ_1,\cR_1}
=&\cK_{n-1}^{\cF_1,\cQ_1,\cR_1}+M_n^{\cQ_1}\eta_n+V_n^{\cQ_1}(\eta_n^2-\omega_n)\\
=&\cK_{n-1}^{\cF_2,\cQ,\cR_2}+M_n^{\cQ}(\eta_n-\eta_{n-1}-(-\eta_{n-1}))\\
&+V_n^{Q}((\eta_n-\eta_{n-1}-(-\eta_{n-1}))^2-\omega_n)\\
=&\cK_n^{\cF_2,\cQ,\cR_2},\\
\cK_n^{\cF_1,\cQ_2,\cR_1}
=&\cK_{n-1}^{\cF_1,\cQ_2,\cR_1}+M_n^{\cQ_1}\eta_n+V_n^{\cQ_1}(\eta_n^2-\omega_n)\\
=&\cK_{n-1}^{\cF_2,\cQ,\cR_3}-M_n^{\cQ}\eta_n+V_n^{\cQ}(\eta_n^2-\omega_n)\\
=&\cK_{n-1}^{\cF_2,\cQ,\cR_3}+M_n^{\cQ}(-\eta_n-\eta_{n-1}-(-\eta_{n-1}))\\
&+V_n^{Q}((-\eta_n-\eta_{n-1}-(-\eta_{n-1}))^2-\omega_n)\\
=&\cK_n^{\cF_2,\cQ,\cR_3}.
\end{align*}
This completes the proof.
\end{proof}

Therefore we are interested in when we can determine the convergence of $\sum_n x_n$ by $\sum_n m_n$ and $\sum_n v_n$,  and when we can not.
From Theorem \ref{th:bound} 
combined with Proposition \ref{prop:bound-ht},
we can easily prove the following relations.

\begin{corollary}\label{cor:bounded-split}
In BFQH Skeptic can force the following events:
\begin{enumerate}
\item 
$\sum_n m_n\textrm{ converges and }\sum_n v_n<\infty\Rightarrow\sum_n x_n\textrm{ converges},$
\item 
$\sum_n m_n\textrm{ does not converge and }\sum_n v_n<\infty\Rightarrow\sum_n x_n\textrm{ does not converge}$,
\item 
$\sum_n m_n\textrm{ converges and }\sum_n v_n=\infty\Rightarrow\sum_n x_n\textrm{ does not converge}.$
\item 
$
\sum_n x_n \textrm{ converges}  \Rightarrow  
\big( \ (\sum_n m_n \textrm{ converges and } \sum_n v_n < \infty ) \textrm{ or } 
$
\begin{equation}
\quad  (\limsup_n \sum_{k=1}^n m_k=\infty, \ 
\liminf_n \sum_{k=1}^n m_k=-\infty  \textrm{ and } 
\sum_n v_n = \infty ) \ \big)
\label{eq:three-series-divergence0}
\end{equation}
\end{enumerate}
\end{corollary}
(i) and (ii) follow from ($\Rightarrow$) of Theorem \ref{th:bound}, 
(iii) follows from ($\Leftarrow$) of Theorem \ref{th:bound} and
(iv) follows from the fact that if Skeptic can force an event $E$ in
the bounded forecasting game, then he can force $E$ in BFQH.

In the classical three-series theorem, the second case
of the right-hand side of \eqref{eq:three-series-divergence0} is 
eliminated by the assumption of independence of the random variables.
In view of Gilat's counter example, it seems that a simple general
statement for the game-theoretic framework can not be given for this case.
However in some special cases, where the behaviors of $\sum_n m_n$ and $\sum_n v_n$ are 
simple, we can give definite statements.  
In Corollary \ref{cor:one-sided BFQH} 
and Corollary  \ref{cor:threepoint}
we discuss such cases.

One simple case is that Reality's move $x_n$ is restricted to be
non-negative.
\begin{corollary} (One-sided BFQH)\ 
\label{cor:one-sided BFQH}
Consider the following special case: $x_n \in [0,C]$ in BFQH.
Then Skeptic can force
\[
\sum_n m_n \text{ converges and } \sum_n v_n \text{ converges }
\ \Leftrightarrow  \ 
\sum_n x_n \text{ converges } 
\]
\end{corollary}
This corollary easily follows because by coherence $m_n \ge 0$ and
and $\sum_n m_n$ is a non-negative series, which
eliminates the second case on the right-hand side of \eqref{eq:xnconvergence}.

We now consider the case that the move space of Reality is restricted to be
a set of three points (trinomial game, cf.\ \cite{nakajima-etal}).  This case will play an essential role in Section \ref{subsec:div}.
Indeed the counter examples to SLLN in 
Section 4.3 of \cite{ShaVov01} and Section 7 of \cite{KumTakTak07}
are constructed as probability distributions on a set of three points.

\begin{corollary}\label{cor:threepoint} (Trinomial Game)\ \ 
We consider the following special case in BFQH:
\[m_n=0,\ v_n\in[0,1],\ x_n\in\{0,\pm1\}.\]
Then Skeptic can force
\[\sum_n v_n<\infty\iff x_n=0\textrm{ for all but finite }n.\]
\end{corollary}

\begin{proof}
This follows from Theorem \ref{th:bound} because for $x_n\in \{0,\pm 1\}$, $\sum_n x_n$ converges if and only if $x_n=0$ for all but finite $n$.
\end{proof}

For the rest of this section we again consider the coherence of BFQH
mentioned just after Theorem \ref{th:bound} in relation to 
Corollary \ref{cor:threepoint}.
The coherence can also be proved by a direct calculation as follows.
If $V_n\ge0$ then Reality announces $x_n=m_n$.
Then $M_n(x_n-m_n)+V_n((x_n-m_n)^2-v_n)\le0$.
Suppose  $V_n<0$.
If $m_n\ge\frac{M_n}{2V_n}$ then Reality announces $x_n=-C$.
Then
\begin{align*}
&M_n(x_n-m_n)+V_n((x_n-m_n)^2-v_n)\\
\le &M_n(-C-m_n)+V_n((-C-m_n)^2-(C^2-m_n^2))\\
=&-M_n(C+m_n)+2V_nm_n(C+m_n)\\
=&(2V_nm_n-M_n)(C+m_n)\le0.
\end{align*}
If $V_n<0$ and $m_n\le\frac{M_n}{2V_n}$ then Reality announces $x_n=C$.
The calculation is the same as above.
This fact is important for our argument in Section \ref{subsec:div}, so we state it as a proposition.
\begin{proposition}\label{th:boundreal}
BFQH remains coherent
even with the restriction $x_n\in\{m_n,\pm C\}$.
\end{proposition}

\section{The rate of convergence of SLLN}
\label{sec:rate}

In this section we consider the rate of convergence of SLLN in the
usual unbounded forecasting game with quadratic hedge.

\begin{quote}
{\sc Unbounded Forecasting}\\
\textbf{Players}: Forecaster, Skeptic, Reality\\
\textbf{Protocol}:\\
\indent $\cK_0=1$.\\
\indent FOR $n=1,2,\ldots$:\\
\indent\indent Forecaster announces $m_n\in\mathbb{R}$ and $v_n > 0$.\\
\indent\indent Skeptic announces $M_n\in\mathbb{R}$ and $V_n\ge0$.\\
\indent\indent Reality announces $x_n\in\mathbb{R}$.\\
\indent\indent $\cK_n:=\cK_{n-1}+M_n(x_n-m_n)+V_n((x_n-m_n)^2-v_n)$.\\
\textbf{Collateral Duties}:
Skeptic must keep $\cK_n$ non-negative.
Reality must keep $\cK_n$ from tending to infinity.
\end{quote}

As in the last section we use the notation
$\ss_n=\sum_{k=1}^n(x_k-m_k)$, $A_n=\sum_{k=1}^n v_k$ and
$A_\infty=\lim_n A_n.$
In many cases we assume that $m_n=0$ for all $n$ without loss of generality.

\subsection{Motivation}\label{subsec:motivation}

The rate of convergence of SLLN for i.i.d.\ case was completely solved by LIL of
Hartman and Wintner \cite{HarWin41}.
We refer to \cite{Pet95}.

\begin{theorem}[Hartman-Wintner's law of the iterated logarithm]
Let $\{X_n\}$ be a sequence of independent identically distributed random variables with zero mean and finite variance $\sigma^2$.
We put $S_n=\sum_{k=1}^n X_k$, $a_n=(2n\log \log n)^{1/2}$.
Then
\[\limsup S_n/a_n=\sigma\mbox{ a.s., }\liminf S_n/a_n=-\sigma\mbox{ a.s.}\]
\end{theorem}

The theorem was a generalization of the case of binary sequences by Khinchin \cite{Khi24}.
If we drop the condition of i.i.d., we need some additional conditions.

\begin{theorem}[Kolmogorov \cite{Kol29}]
\label{thm:kolmogorov29}
Let $\{X_n\}$ be a sequence of independent random variables with zero means and finite variances.
Put $\sigma^2_n=\textrm{Var}X_n$ and $B_n=\sum_{k=1}^n \sigma^2_k$.
Suppose $B_n\to\infty$.
Suppose also that there exists a sequence of positive constants $\{M_n\}$ such that
\[M_n=o\left(\left(\frac{B_n}{\log\log B_n}\right)^{1/2}\right)\]
and
\[|X_n|\le M_n\textrm{ a.s.}\]
Then
\[\limsup\frac{S_n}{(2B_n\log\log B_n)^{1/2}}=1\textrm{ a.s.}\]
\end{theorem}

Some other condition for LIL than in  Theorem \ref{thm:kolmogorov29}
is given in \cite{tomkins-1978} and LIL for martingales are discussed for example in 
\cite{stout-1970} and \cite{fisher-1992}.
For game-theoretic LIL see Chapter 5 of 
\cite{ShaVov01} and a recent paper by Takazawa \cite{takazawa-aism}.

In game-theoretic probability we can not assume that the sequence $\{x_n\}$ announced by Reality is i.i.d.
Furthermore we can not assume the existence of $M_n$ either such that $|x_n|\le M_n$.
From now on we consider the rate of convergence of SLLN in game-theoretic probability.
In the view point of measure-theoretic probability it is the rate of convergence of SLLN in a non-identical case.

\subsection{Results on convergences in unbounded forecasting}\label{subsec:conv}

Here we derive several results on the rate of convergence of
Kronecker's lemma and hence the strong law of large numbers.

\begin{lemma}\label{lem:vconv}  In the unbounded forecasting 
Skeptic can force
\[\sum_n v_n<\infty\Rightarrow\sum_{k=1}^n (x_k-m_k)\textrm{ converges.}\]
Skeptic can not force 
\[
\sum_{k=1}^n (x_k-m_k)\textrm{ converges} \Rightarrow \sum_n v_n<\infty.
\]
\end{lemma}

\begin{proof}
The proof of the first statement is exactly the same as the proof of 
($\Rightarrow$) in Theorem \ref{th:bound}. For the second statement,
consider Reality's strategy $x_n=m_n$, $\forall n$.
Then $\cK_n = \cK_0 - \sum_{k=1}^n V_k v_k$ and clearly
Skeptic has no control over $v_n$'s and hence
can not achieve $\sum_n v_n < \infty$.
\end{proof}

\begin{theorem}\label{th:conv}
Let $g$ be a positive increasing function. In the unbounded forecasting
Skeptic can force
\begin{align}\label{eq:conv}
\sum_n\frac{v_n}{g(A_n)}<\infty\Rightarrow \sum_{k\le n}\frac{x_k-m_k}{\sqrt{g(A_k)}}\textrm{ converges.}
\end{align}
\end{theorem}

\begin{proof}
We assume that $m_n=0$ for all $n$ without loss of generality.
We consider the capital process
\[W_n=\sum_{k\le n}\frac{x_k}{\sqrt{g(A_k)}}.\]
The compensator of $W^2$ is 
\[B_n=\sum_{k=1}^n \frac{v_k}{g(A_k)}.\]
If $B_\infty < \infty$, then 
$W_n$ converges by Lemma 4.6 and Lemma 4.7 in \cite{ShaVov01}.
\end{proof}

\begin{corollary}
\label{cor:convergent-g-integral}
Let $g$ be a positive increasing function on $[0,\infty)$ with $g(\infty)=\infty$.
Skeptic can force
\begin{align}\label{eq:vinfconv}
\sum_n v_n=\infty\textrm{ and }\sum_n\frac{v_n}{g(A_n)}<\infty
\Rightarrow \frac{\sum_{k\le n}(x_k-m_k)}{\sqrt{g(A_n)}}\to0.
\end{align}
\end{corollary}

\begin{proof}
This follows easily from Theorem \ref{th:conv} and Kronecker's lemma.
\end{proof}

In some cases we can drop $\sum_n v_n/g(A_n)<\infty$ from (\ref{eq:vinfconv}).
The following is an example.

\begin{corollary}\label{cor:conv}
Let $g$ be a positive increasing function such that
\[\int_0^\infty \frac{1}{g(x)}dx<\infty.\]
Then Skeptic can force
\[A_\infty=\infty\Rightarrow \frac{\sum_{k\le n}(x_k-m_k)}{\sqrt{g(A_n)}}\to0.\]
\end{corollary}

\begin{proof}
It suffices to show that $\sum_n v_n/g(A_n)<\infty$ when $A_\infty=\infty$.
This holds because
\begin{align*}
\int_0^\infty\frac{1}{g(x)}dx
=\sum_{n=1}^\infty\int_{A_{n-1}}^{A_n}\frac{1}{g(x)}dx
\ge\sum_{n=1}^\infty\frac{v_n}{g(A_n)}.
\end{align*}
\end{proof}

\begin{example}
\label{ex:iterated-log}
We write $\ln^i$ to mean the function such that $\ln^i(x)=\ln(\ln^{i-1}(x))$ and $\ln^0(x)=x$ defined recursively.
Let
\[g_i(x):=(\prod_{j=0}^i\ln^j)\times\ln^i(x).\]
In other words,
\[g_0(x)=x^2,\ g_1(x)=x(\ln x)^2,\ g_2(x)=x\ln x(\ln\ln x)^2.\]
Then
\[\int \frac{1}{g_i(x)}dx=-\frac{1}{\ln^i x}.\]
Hence Skeptic can force
\[A_\infty=\infty\Rightarrow \frac{\sum_{k\le n}(x_k-m_k)}{\sqrt{g_i(A_n)}}\to0\]
for all $i$.

For example, consider the special case $v_n \equiv v$. Then $A_\infty=\infty$ is automatic and
with $g_2(x)$ above we have
$\sum_{k\le n} (x_k-m_k)/(\sqrt{n\ln n} \ln\ln n) \rightarrow 0$.
\end{example}

Note that it does not follow $\sum_{k\le n}(x_k-m_k)/\sqrt{n \ln n \ln\ln n}\rightarrow 0$.
This is because $\sum_n \frac{1}{n\ln n \ln\ln n}=\infty$.
However this does not mean that Theorem \ref{th:conv} is a weaker result compared to LIL.
Recall that we are considering a non-identical case in the measure-theoretic point of view.
In fact Theorem \ref{th:conv} is strict as we will see in the next subsection.

We now apply Corollary \ref{cor:convergent-g-integral}
to the coin-tossing game.  We first show that
the coin-tossing game is a special case of BFQH.
Restrict $x_n \in \{0,1\}$, $m_n=p_n \in [0,1]$ and $v_n = p_n (1-p_n)$ in
BFQH.
For $M, V\in {\mathbb R}$, $p\in [0,1]$, consider
\[
M (x-p) + V ((x-p)^2 - p(1-p))
\]
For  both $x=0$ and $x=1$  we have $(x-p)^2 - p(1-p)=(x-p) (1-2p)$.  Therefore
\begin{align*}
M (x-p) + V ((x-p)^2 - p(1-p)) &= M(x-p) +  V(x-p) (1-2p)\\
& = (M+V(1-2p)) (x-p).
\end{align*}
By considering Skeptic's move $M_n + V_n(1-2p_n)$, we see that
BFQH for $x_n \in \{0,1\}$ reduces to the coin-tossing game.
Therefore from Corollary \ref{cor:conv} applied to the coin-tossing game,
we have the following result. 

\begin{proposition}
\label{prop:levy-ii}
Let $g$ be a positive increasing function such that
$\int_0^\infty \frac{1}{g(x)}dx<\infty$.
Let $\bar A_n = \sum_{k=1}^n p_k$.
Then in the coin-tossing game Skeptic can force
\[\bar A_\infty=\infty\Rightarrow \frac{\sum_{k\le n}(x_k-p_k)}{\sqrt{g(\bar A_n)}}\to0.\]
\end{proposition}

\begin{proof}
Let $A_n = \sum_{k=1}^n p_k(1-p_k)$. Then $A_n\le \bar A_n$.
If $A_\infty=\bar A_\infty =\infty$, then by Corollary 
\ref{cor:conv} Skeptic can force 
\[
\frac{\sum_{k\le n}(x_k-p_k)}{\sqrt{g(\bar A_n)}}
=\frac{\sum_{k\le n}(x_k-p_k)}{\sqrt{g(A_n)}} \sqrt{\frac{g(A_n)}{g(\bar A_n)}}
\to 0.
\]
On the other hand, if $A_\infty < \bar A_\infty=
\infty$, then Skeptic can force the convergence of $\sum_n (x_n-p_n)$.
However in this case $\sum_{k\le n}(x_k-p_k)/\sqrt{g(\bar A_n)}\to0$, 
because $g(\infty)=\infty$.
\end{proof}

Note that Proposition \ref{prop:levy-ii}  strengthens  (ii) of
Remark \ref{pro:levy}, which only states
\[
\bar A_\infty=\infty\Rightarrow \frac{\sum_{k\le n}(x_k-p_k)}{\bar A_n}
= \frac{\sum_{k\le n}x_k}{\sum_{k\le n} p_k } - 1
\to0.
\]

\subsection{Results on divergence in unbounded forecasting}\label{subsec:div}

As the converse to the convergence result in the previous subsection,
we will prove that Skeptic can not force  the convergence on 
the right-hand side of (\ref{eq:conv}) when $\sum_n v_n/g(A_n)=\infty$.
The novelty in our approach is that in order to prove this fact we
use deterministic strategies of Reality.  

We formulate Reality's strategies and 
introduce the notion of {\it compliance} of Reality with an event.
We propose to use the term ``compliance'' for Reality's strategies and
reserve the word ``forcing'' to Skeptic for clarity of our arguments.
Skeptic's forcing an event $E$ means that
Reality {\it has to} move so that $E$ happens
in order that Skeptic's capital stays bounded.
In contrast Reality's forcing an event $E$ means intuitively that
Reality {\it can} move so that $E$ happens and Skeptic's capital stays bounded
irrespective of the moves of Forecaster and Skeptic.
In other words Reality conforms to $E$ and adapts as requested.

For notational simplicity, as in the multi-dimensional bounded forecasting
game, write $\mu_n=(m_n, v_n)$, $\fm_n = (M_n, V_n)$.
Consider Reality's strategy $\cR=\{\cR_n\}_{n\ge 1}$ which 
determines Reality's move $x_n$ based on the moves $\mu_k,\fm_k$, $k\le n$,
of Forecaster and Skeptic: 
\[
x_n = \cR_n(\mu_1, \fm_1, \dots, \mu_n, \fm_n), \quad n\ge 1.
\]

\begin{definition}
\label{def:reality-stragety}
By a strategy $\cR$ Reality {\em complies} with the event $E\subset \Xi$, 
if  
\begin{itemize}
\item[(i)] irrespective of the moves $\mu_n, \fm_n$, $n\ge 1$, of Forecaster 
and Skeptic, both observing their collateral duties, it holds that
\[
(\mu_1,\cR_1(\mu_1,\fm_1), \mu_2,\cR_2(\mu_1,\fm_1, \mu_2, \fm_2), \dots ) \in E,
\]
and 
\item[(ii)]
$\sup_n \cK_n < \infty$.
\end{itemize}
We say that by $\cR$ Reality {\em strongly} complies with $E$ 
if the supremum in (ii) is uniformly bounded from above by $1=\cK_0$, i.e., 
$\cK_n \le 1$ %
irrespective of the moves of Forecaster 
and Skeptic, both observing their collateral duties.

We simply say that Reality 
(strongly) complies 
with the event $E\subset \Xi$ if there exists a strategy $\cR$
such that by $\cR$ Reality (strongly) complies with $E$.
\end{definition}

Concerning the notion of compliance we prove the following fundamental proposition.

\begin{proposition}\label{lem:SRforce}
In the unbounded forecasting, 
if Skeptic can force an event $E$, then Reality strongly complies with $E$.
\end{proposition}

\begin{remark}
\label{rem:SRforce}
As seen in the proof of this proposition below, the statement
holds not only for the unbounded forecasting, but for more
general protocols of game-theoretic probability, including the trinomial game.
\end{remark}

The idea of the proof is as follows.
First of all Reality needs to prevent the capital from tending to infinity.
By coherence this is possible for Reality.
Next the path must be in $E$.
Since Skeptic can force $E$, Skeptic has the strategy such that if the path is not in $E$, then the capital goes to infinity.
It follows that if the capital does not tend to infinity, then the path is in $E$.
Then all Reality has to do is to prevent the capital of the strategy from tending to infinity.
Again by coherence this is possible for Reality.
Can Reality prevent the capitals of two strategies from tending to infinity?
It is possible by considering the strategy that is the average of two strategies.
In other words Reality's strategy can be constructed by considering a single sufficiently powerful strategy of Skeptic.
Furthermore Reality's strategy can be deterministic.
To make the strategy ``strongly'' comply we need a more precise argument as is in the proof.

Such an argument is commonly used in algorithmic randomness.
Especially some examples of random sets are sets on which a single sufficiently powerful (super)martingale fails.
See \cite{Nie09,Gac85,MerMih02,NST05,Miyabe_extvan}.
One way to obtain a set on which the (super)martingale fails is to choose the leftmost non-ascending path in binary sequences.
This choice corresponds to the coherence in game-theoretic probability.
Although the constructed random sets may not be computable in general, it can be constructed deterministically by the (super)martingale.

We set up some more notation for clarity.
When the moves of all the players are individually specified 
we write Skeptic's capital as 
\[
\cK_n[(\mu_k,\fm_k,x_k)_{k=1}^n)],
\quad (\cK_0=1).
\]
In this notation 
Skeptic's capital under a strategy $\cP$ is written as
\[
1+\cK_n^\cP = \cK_n [(\mu_k,\cP_k,x_k)_{k=1}^n].
\]
We now give a proof of Proposition \ref{lem:SRforce}.

\begin{proof}
Since Skeptic can force $E$, there exists Skeptic's strategy $\cP$ such that
\begin{align}\label{eq:forceA}
\limsup_n\cK_n^\cP<\infty\Rightarrow (\mu_n,x_n)_{n=1}^\infty \in E.
\end{align}
First we give a strategy $\cR$ of Reality 
such that %
$\cK_n$ 
is uniformly bounded from above by $1+\epsilon$, where
$\epsilon>0$ is arbitrarily fixed.

Consider Reality's move at the first round $n=1$ after Forecaster's move
$\mu_1=(m_1, v_1)$ and Reality's move $\fm_1=(M_1,V_1)$ were announced.  Write 
$\cP_1(\mu_1)=\cP_1((m_1,v_1))=(M_1^\cP, V_1^{\cP})$,
which is the move of the strategy $\cP$ at the first round.  Let $\alpha=1/(1+\epsilon)$ and let
\begin{equation*}
\tilde \fm_1=(\tilde M_1, \tilde V_1)= (1-\alpha)(M_1^\cP, V_1^{\cP}) + \alpha (M_1, V_1).
\end{equation*}
Because of coherence, Reality can (deterministically) 
choose $x_1$ such that
\begin{align*}
\cK_1[\mu_1, \tilde \fm_1, x_1]
&=
(1-\alpha) \cK_1[\mu_1, \cP_1(\mu_1), x_1] 
 +\alpha \cK_1[\mu_1, \fm_1,  x_1] \\
&\le\cK_0=1.
\end{align*}
Since both 
$\cK_1[\mu_1, \cP_1(\mu_1), x_1], \cK_1[\mu_1, \fm_1,  x_1]$
are non-negative, it follows that 
\[
\cK_1[\mu_1, \cP_1(\mu_1),  x_1]\le \frac{1}{1-\alpha}=\frac{1+\epsilon}{\epsilon}, \quad 
\cK_1[\mu_1, \fm_1,  x_1]\le \frac{1}{\alpha}=1+\epsilon.
\]

We now make an inductive argument.  Suppose that
Reality could deterministically choose $x_1,\dots, x_{n-1}$ such that
\[
\cK_{n-1}[(\mu_k, \cP_k,  x_k)_{k=1}^{n-1}] \le \frac{1+\epsilon}{\epsilon}, \ \ 
\cK_{n-1}[(\mu_k, \fm_k,  x_k)_{k=1}^{n-1}] \le 1+\epsilon.  %
\]
As in the first round define
\begin{equation*}
\tilde \fm_n= (\tilde M_n, \tilde V_n)= (1-\alpha)(M_n^\cP, V_n^{\cP}) + \alpha (M_n, V_n),
\end{equation*}
where $(M_n, V_n)=\fm_n$ is the actual move announced by Skeptic
and $(M_n^\cP, V_n^{\cP})$ is the move of strategy $\cP$.
By coherence, Reality can now choose $x_n$ such that
\[
\tilde M_n (x_n - m_n) + \tilde V_n ( (x_n - m_n)^2 - v_n) \le 0.
\]
Then 
\begin{align}
\cK_n[(\mu_k, \tilde \fm_k, x_k)_{k=1}^n]
&=(1-\alpha) \cK_n[(\mu_k, \cP_k, x_k)_{k=1}^n] 
 +\alpha\cK_n[(\mu_k, \fm_k,  x_k)_{k=1}^n] \nonumber \\
&\le\cK_{n-1}[(\mu_k, \tilde \fm_k, x_k)_{k=1}^{n-1}] \nonumber \\
&\le 1.
\label{eq:almost-strong-compliance0}
\end{align}
Thus as in the first round 
\begin{equation}
\label{eq:both-bounded}
\cK_n[(\mu_k, \cP_k, x_k)_{k=1}^n] \le \frac{1+\epsilon}{\epsilon}, \  \ 
\cK_n[(\mu_k, \fm_k,  x_k)_{k=1}^n]\le 1+\epsilon.
\end{equation}
By \eqref{eq:forceA} and the first term of 
\eqref{eq:both-bounded}, (i) of Definition
\ref{def:reality-stragety} is satisfied.
By the second term of 
\eqref{eq:both-bounded}, 
$\cK_n[(\mu_k, \fm_k,  x_k)_{k=1}^n]$ is uniformly bounded from above by
$1+\epsilon$.

It remains to show that we can let $\epsilon=0$.
We argue as follows. 
By coherence, Reality can always choose $x_n$ such that 
$M_n (x_n - m_n) + V_n ( (x_n - m_n)^2 - v_n) \le 0.$
In the unbounded forecasting, unless $(M_n, V_n)=(0,0)$, 
Reality can choose $x_n$ such that this inequality is strict.
Reality will look for the first time $n=n_0$  such that $(M_n, V_n)\neq (0,0)$.
At round $n_0$ Reality chooses $x_{n_0}$ such that
$\cK_{n_0}[(\mu_k, \fm_k,  x_k)_{k=1}^{n_0}] < 1$.
For $n=1,\dots,n_0-1$, Reality chooses  $x_n$ such that
$\cK_n^\cP \le 0$. 
Now define $\alpha= \cK_{n_0}[(\mu_k, \fm_k,  x_k)_{k=1}^{n_0}]$ and after the round $n_0$
Reality follows the strategy ensuring $\cK_n \le 1$, $n>n_0$, 
as in 
\eqref{eq:almost-strong-compliance0}.
On the other hand, if there is no such $n_0$, then Reality keeps 
choosing  $x_n$ such that $\cK_n^\cP \le 0$. In this case
$\cK_{n}[(\mu_k, \fm_k,  x_k)_{k=1}^{n}]=1$ for all $n$ %
and also (i) of Definition \ref{def:reality-stragety} is satisfied by \eqref{eq:forceA}.
\end{proof}

We now state the following theorem.

\begin{theorem}\label{th:notconv}
Let $g:\mathbb{R}\to\mathbb{R}$ be a positive increasing function.
Then in the unbounded forecasting Reality strongly complies with the event
\begin{align}\label{eq:notconv}
\sum_n\frac{v_n}{g(A_n)}=\infty\Rightarrow \frac{\sum_{k\le n}(x_k-m_k)}{\sqrt{g(A_n)}}\textrm{ does not converge.}
\end{align}
\end{theorem}

As an immediate consequence 
of this theorem we have the following corollary.
\begin{corollary}
\label{cor:notconv}
Let $g:\mathbb{R}\to\mathbb{R}$ be a positive increasing function.
Let $E_1$ be any event depending only on $v_1, v_2, \dots$, such that
$E_1 \cap \{\sum_n\frac{v_n}{g(A_n)}= \infty\} \neq\emptyset$.
In the unbounded forecasting Skeptic can not force
\[
E_1 \Rightarrow \frac{\sum_{k\le n}(x_k-m_k)}{\sqrt{g(A_n)}}\textrm{ converges}.
\]
\end{corollary}

For proving Theorem \ref{th:notconv} we prove two technical lemmas.

\begin{lemma}\label{lem:convbound}
Let $\{y_n\}$ be a sequence of reals and let $\{g_n\}$ be 
a non-decreasing sequence of positive reals. %
If $(\sum_{k\le n}y_k)/g_n$ converges to $d$, 
then $|y_n/g_n|\le|d|+1$ for all but finite $n$.
\end{lemma}
\begin{proof}
First note that
\begin{align*}
\left|\frac{y_n}{g_n}\right|
\le&\left|\frac{y_n}{g_n}-
d \left(1-\frac{g_{n-1}}{g_n}\right)\right|+
\left|d \left(1-\frac{g_{n-1}}{g_n}\right)\right|
\le |d| + \left|\frac{y_n}{g_n}-
d \left(1-\frac{g_{n-1}}{g_n}\right)\right|.
\end{align*}
Therefore it suffices to show  that 
\[\left|\frac{y_n}{g_n}- d \left(1-\frac{g_{n-1}}{g_n}\right)\right|\le1\]
for all sufficiently large $n$.

Let $\epsilon$ be such that $0<\epsilon\le\frac{1}{3}$.
Then there exists $N$ such that 
\begin{align*}
n>N\Rightarrow \left|\frac{\sum_{k=1}^{n}y_k}{g_n}-d\right|<\epsilon.
\end{align*}
It follows that, for all $n-1>N$,
\begin{align*}
&\left|\frac{y_n}{g_n}-d \cdot\left(1-\frac{g_{n-1}}{g_n}\right)\right|\\
\le&\left|\frac{y_n}{g_n}-\frac{\sum_{k=1}^{n-1} y_k}{g_{n-1}}\cdot\left(1-\frac{g_{n-1}}{g_n}\right)\right|
+\left|\frac{\sum_{k=1}^{n-1} y_k}{g_{n-1}}-
d \right|\cdot\left(1-\frac{g_{n-1}}{g_n}\right)\\
<&\left|\frac{\sum_{k=1}^n y_k}{g_n}-\frac{\sum_{k=1}^{n-1} y_k}{g_{n-1}}\right|+\epsilon
\le\left|\frac{\sum_{k=1}^n y_k}{g_n}-d\right|+\left|\frac{\sum_{k=1}^{n-1} y_k}{g_{n-1}}-d\right|+\epsilon<3\epsilon .
\end{align*}
\end{proof}

\begin{lemma}\label{lem:ser}
Let $\{a_n\}$ be a sequence of positive reals.
Then there exists a sequence $\{\epsilon_n\}$ of positive reals such that
\begin{enumerate}
\item $\epsilon_n$ is determined only by $a_1,\cdots,a_n$,
\item $\epsilon_n a_n\le1$,
\item $\sum_n a_n=\infty$ implies $\sum_n \epsilon_n a_n=\infty$ and $\epsilon_n\to0$.
\end{enumerate}
\end{lemma}

\begin{proof}
We define $\epsilon_n$ as follows.
\begin{enumerate}
\item[(P1)] Let $n=b=c=1$.
\item[(P2)] If $2^{-b}a_n\ge1$ then let $\epsilon_n=1/a_n$, otherwise $\epsilon_n=2^{-b}$, i.e.\ 
let $\epsilon_n = \min(1/a^n, 2^{-b})$.
\item[(P3)] If $\sum_{k=c}^n \epsilon_k a_k\ge1$ then let $b=b+1,$ and $c=n+1$.
\item[(P4)] Let $n=n+1$ and goto (P2).
\end{enumerate}

It is clear that (i) and (ii) are satisfied.
We shall prove that (iii) is satisfied.
Suppose that $\sum_n a_n=\infty$.

We claim that $\sum_{k=c}^n\epsilon_k a_k\ge1$ in (P3) for infinitely many times.
Otherwise, $b$ and $c$ do not change from some point.
Also, in view of 
\[
\sum_{k=c}^n \epsilon _k a_k < 1 \ \Rightarrow \ \epsilon_n a_n < 1 \ \Rightarrow \epsilon_n = e^{-b},
\]
$\epsilon_n=2^{-b}$ for all but finite $n$.
It follows that $\sum_{k\ge c}2^{-b}a_k<1$.
This contradicts to the fact that $\sum_n a_n=\infty$.

Then there exists an increasing sequence $\{c_i\}$ such that $\sum_{k=c_i+1}^{c_{i+1}}\epsilon_k a_k\ge1$.
It follows that $\sum_n\epsilon_n a_n=\infty$.

Since $b$ goes to infinity and $\epsilon_n\le2^{-b}$ in (P2), we have $\epsilon_n\to0$.
\end{proof}

We are now ready to give  a proof of Theorem \ref{th:notconv}.
Before starting a formal proof, we discuss the idea of the proof.
By Proposition \ref{lem:SRforce} we wish Skeptic could force the divergence.
As we saw in the Lemma \ref{lem:vconv}, however, it is not easy for Skeptic to force the divergence in the unbounded forecasting.
So we use the bounded forecasting, more precisely, the game restricted to three points as in Proposition \ref{th:boundreal}.

It is suffices to show that Skeptic can force (\ref{eq:notconv}) in the restricted protocol by the following reason.
Consider the following three statements.
\begin{enumerate}
\item Skeptic can force the event in the restricted protocol.
\item Reality strongly complies with the event in the restricted protocol.
\item Reality strongly complies with the event in the unbounded forecasting.
\end{enumerate}
The implication (ii)$\Rightarrow$(iii) holds by the following simple argument.
Suppose that Reality strongly complies with it in the restricted protocol.
Then Reality can use the same strategy in the unrestricted protocol or in the unbounded forecasting.

The implication (i)$\Rightarrow$(ii) follows from Proposition \ref{lem:SRforce}.
Note that the restricted protocol Skeptic can use a strategy with $V_n\le0$ 
in order to force (\ref{eq:notconv}).
In the unrestricted protocol this is not allowed.
It is Reality who refers this Skeptic's strategy to make her strategy.

To prove (i), we use the trinomial game 
(cf.\ Corollary \ref{cor:threepoint} and Proposition \ref{th:boundreal}).
In other words we construct a reduction from the restricted protocol to the trinomial game.
Then since Skeptic can force the divergence in the trinomial game, Skeptic can also force the divergence in the restricted protocol.
This argument of the reduction is often seen in computability theory \cite{Odi90,Odi99}, complexity theory \cite{AroBar09}
and algorithmic randomness \cite{Dow10,Nie09}.

We now give a formal proof.

\begin{proof}[Proof of Theorem \ref{th:notconv}]
For notational simplicity we write
$y_n=x_n-m_n$, $g_n=\sqrt{g(A_n)}$ and $z_n=y_n/g_n$.
Let $\epsilon_n\ge0$ be such that $\epsilon_n^2$ is a sequence of Lemma \ref{lem:ser} for $a_n=\frac{v_n}{g(A_n)}$.
We restrict Reality's moves as $\epsilon_n z_n\in\{0,\pm1\}$.
The capital process is
\begin{align*}
\cK_n
&=\cK_{n-1}+M_n(x_n-m_n)+V_n((x_n-m_n)^2-v_n)\\
&=\cK_{n-1}+M_n \frac{g_n}{\epsilon_n} \epsilon_n z_n+V_n \frac{g(A_n)}{\epsilon_n^2}\left((\epsilon_n z_n)^2-\epsilon_n^2\frac{v_n}{g(A_n)}\right)\\
&=\cK_{n-1}+M'_n(x'_n-m'_n)+V'_n((x'_n-m'_n)^2-v'_n).
\end{align*}
Then $\cK_n$ can be seen as a capital process in Corollary \ref{cor:threepoint} where Forecaster announces $m'_n=0$ and $v'_n=\epsilon_n^2\frac{v_n}{g(A_n)}$, 
Skeptic announces $M'_n=M_n \frac{g_n}{\epsilon_n}$ and $V'_n=V_n \frac{g(A_n)}{\epsilon_n^2}$ and Reality announces $x'_n=\epsilon_n z_n\in\{0,\pm1\}$.
Note that $v'_n\in[0,1]$ because of the definition of $\epsilon_n$.
Since $|\epsilon_n z_n|$ is bounded, $V_n$ can also be negative in the restricted protocol.
Hence, by Corollary \ref{cor:threepoint}, Skeptic can force that 
\begin{align*}
\sum_n\epsilon_n^2\frac{v_n}{g(A_n)}=\infty
\Rightarrow \epsilon_n z_n=\pm1\textrm{ for infinitely many }n
\end{align*}
and Reality can choose her move such that 
$\cK_n\le \cK_{n-1}\le \dots \le \cK_0=1$ for all $n$.

By the definition of $\epsilon_n$ we have $\sum_n\frac{v_n}{g(A_n)}=\infty\Rightarrow\sum_n\epsilon_n^2\frac{v_n}{g(A_n)}=\infty$.
Furthermore If $\epsilon_n|z_n|=1$ for infinitely many $n$ then
$|z_n|=1/\epsilon_n$ for infinitely many $n$ and $\sup_n|z_n|=\infty$ by $\epsilon_n^2\to0$ and $\epsilon_n\to0$.
It follows that $(\sum_{k\le n}y_k)/g_n$ does not converge by Lemma \ref{lem:convbound}.
\end{proof}

We can summarize the results of this section as follows.
Consider a game in which Skeptic wins when $\frac{\sum_{k\le n}(x_k-m_k)}{\sqrt{g(A_n)}}$ converges.
If both players do the best, then the winner depends only on whether $\sum_n\frac{v_n}{g(A_n)}<\infty$ or not.
In other words in order to make it converge we need a function $g$ such that $\sum_n\frac{v_n}{g(A_n)}<\infty$.
As we have already seen in Example \ref{ex:iterated-log} this $g$ grows faster than $n\log\log n$ that appears in LIL.

A measure-theoretic interpretation is as follows.
Let $m_n$ and $v_n$ be the conditional mean and the conditional variance.
Under any probability measure, 
if $\sum_n\frac{v_n}{g(A_n)}<\infty$, then 
$\sum_{k\le n}\frac{x_k-m_k}{\sqrt{g(A_k)}}$ converges a.s.
Now consider the case $\sum_n\frac{v_n}{g(A_n)}=\infty$.
There exists a  probability measure, such that 
if $\sum_n\frac{v_n}{g(A_n)}=\infty$
then $\sum_{k\le n}\frac{x_k-m_k}{\sqrt{g(A_k)}}$ does not converge a.s.

\section{Rate of convergence of SLLN under a general hedge}
\label{sec:general-hedge}

Theorem \ref{th:conv} can be seen as an extension of SLLN.
In contrast Kumon, Takemura and Takeuchi \cite{KumTakTak07} proved another extension.
Furthermore both extensions have similar forms.
To clarify a relation between these two extensions we consider a protocol
with a general hedge.

\begin{quote}
{\sc Unbounded Forecasting with General Hedge} (UFGH)\\
\textbf{Parameters}: A single hedge $h:\mathbb{R}\to\mathbb{R}$\\
\textbf{Players}: Forecaster, Skeptic, Reality\\
\textbf{Protocol}:\\
\indent $\cK_0=1$.\\
\indent FOR $n=1,2,\ldots$:\\
\indent\indent Forecaster announces $m_n\in\mathbb{R}$ and $v_n > 0$.\\
\indent\indent Skeptic announces $M_n\in\mathbb{R}$ and $V_n\ge0$.\\
\indent\indent Reality announces $x_n\in\mathbb{R}$.\\
\indent\indent $\cK_n:=\cK_{n-1}+M_n(x_n-m_n)+V_n(h(x_n-m_n)-v_n)$.\\
\textbf{Collateral Duties}:
Skeptic must keep $\cK_n$ non-negative.
Reality must keep $\cK_n$ from tending to infinity.
\end{quote}

We assume a few conditions for $h:\mathbb{R}\to\mathbb{R}$.
\begin{itemize}
\item[(A0)] $h(x)=h(|x|)\ge0$.
\item[(A1)] For some $c>0$, $h(x)/x$ is monotone increasing for $x>c$.
\item[(A2)] For some $c>0$, $h(x)/x^2$ is monotone decreasing for $x>c$.
\end{itemize}

In (A2) we are considering hedges weaker than the quadratic hedge, which
is our main interest in this section.  The implications
of hedges stronger than the quadratic hedge to the rate of SLLN
are investigated by recent papers of Takazawa (\cite{takazawa-stochastics}, \cite{takazawa-aism}).

\subsection{The case of convergence}\label{subsec:combconv}

The following theorem is a generalization of Theorem 
\ref{th:conv} to the general hedge.

\begin{theorem}\label{th:comb}
Suppose that $h$ satisfies (A0)-(A2) and that $g$ is a positive increasing function.
Then in UFGH Skeptic can force
\begin{align*}
\sum_n \frac{v_n}{g(A_n)}<\infty\Rightarrow \sum_{k=1}^n\frac{x_k-m_k}{h^{-1}\circ g(A_k)}\textrm{ converges.}
\end{align*}
\end{theorem}

\begin{corollary} In UFGH Skeptic can force
\begin{align*}%
\sum_n v_n<\infty\Rightarrow \sum_{k=1}^n(x_k-m_k)\textrm{ converges.}
\end{align*}
and
\begin{align*}%
\sum_n v_n=\infty\textrm{ and }\sum_n \frac{v_n}{g(A_n)}<\infty\Rightarrow \frac{\sum_{k=1}^n(x_k-m_k)}{h^{-1}\circ g(A_n)}\to0.
\end{align*}
\end{corollary}

This corollary follows from the theorem and Kronecker's lemma.

Notice that Theorem \ref{th:comb} with $h(x)=x^2$ implies Theorem \ref{th:conv}.
In contrast Theorem \ref{th:comb} with $g(x)=h(x/\nu),\ m_n=0$ and $v_n=\nu$ implies Theorem 3.1 of \cite{KumTakTak07}.

The rest of this section is devoted to a proof of Theorem \ref{th:comb}.
The proof is just a straightforward combination of that of Theorem \ref{th:conv} and that of Theorem 3.1 of \cite{KumTakTak07}.
Note that, without loss of generality, we assume that $m_n=0$ for all $n$.

\begin{remark}\label{rem:assh}
We can assume that $c$ in (A1) and (A2) are 0 and $h(x)=x^2$ for $|x|\le1$ by the following reason.
Let $c$ be such that $h(x)/x$ is monotone increasing and $h(x)/x^2$ is monotone decreasing for $x>c$.
Let
\[h_0(x)=\begin{cases}\frac{h(c)}{c}x&\textrm{ if }|x|\le c\\h(x)&\textrm{ otherwise.}\end{cases}\]
Then $h_0$ satisfies (A1) and (A2) for $x\ge0$.
It also follows that $h^{-1}(y)\ge0$ is defined for $y\ge0$.
\end{remark}

By this remark we assume that $c$ in (A1) and (A2) are 0.
Let 
\[h_1(x)=\begin{cases}x^2&\textrm{ if }|x|\le 1\\ h(h^{-1}(1)\cdot x)&\textrm{ otherwise.}\end{cases}\]
Then $h_1$ satisfies (A0)-(A2).
Furthermore, for large $x$, $h_1(x)=h(h^{-1}(1)\cdot x)=y$ and $h^{-1}(y)=h^{-1}(1)\cdot x=h^{-1}(1)\cdot h_1^{-1}(y)$.
Then the convergence does not depend on the use of $h$ and $h_1$.
We also have $h(x)/x^2\le1$.

\begin{lemma}\label{lem:hgconv}
In UFGH Skeptic can force
\begin{align*}
\sum_n \frac{v_n}{g(A_n)}<\infty\Rightarrow \sum_n \frac{h(x_n)}{g(A_n)}<\infty.
\end{align*}
\end{lemma}

\begin{proof}
Consider the strategies such that
\[\cK_0=D,\ M_n=0,\ V_n=\frac{1}{g(A_n)}\]
as long as $\cK_n$ is non-negative for all $x_n$ where $D$ is a natural number.
If $\sum_n v_n/g(A_n)<D$ then the capital process $\cK_n$ is
\begin{align*}
\cK_n
=D+\sum_{k=1}^n \frac{h(x_k)-v_k}{g(A_k)}.
\end{align*}
If $\sum_n v_n/g(A_n) < \infty$ and $\sum_n h(x_n)/g(A_n)=\infty$, then
$\lim_n \cK_n =\infty$.  This proves the lemma.
\end{proof}

\begin{lemma}\label{lem:finite}
Under the conditions (A0),   In UFGH
Skeptic can force
\begin{align*}
\sum_n \frac{v_n}{g(A_n)}<\infty\Rightarrow |x_n|\ge h^{-1}\circ g(A_n)
\end{align*}
for only finitely many $n$.
\end{lemma}

\begin{proof}
Notice that
\[|x_n|\ge h^{-1}\circ g(A_n)\]
implies
\[\frac{h(x_n)}{g(A_n)}\ge1.\]
Then the result follows from Lemma \ref{lem:hgconv}.
\end{proof}

\begin{lemma}
Under the conditions (A0) and (A2), In UFGH Skeptic can force
\begin{align*}
\sum_n \frac{v_n}{g(A_n)}<\infty\Rightarrow \sum_n\frac{x_n^2}{(h^{-1}\circ g(A_n))^2}<\infty.
\end{align*}
\end{lemma}

\begin{proof}
By Lemma \ref{lem:finite} we only consider the case $x_n\le h^{-1}\circ g(A_n)$.
By (A2), $h(x)/x^2$ is monotone decreasing for $x>0$.
Hence
\[\frac{h(x_n)}{x_n^2}\ge\frac{h\circ h^{-1}\circ g(A_n)}{(h^{-1}\circ g(A_n))^2}.\]
It follows that
\[\frac{x_n^2}{(h^{-1}\circ g(A_n))^2}\le\frac{h(x_n)}{g(A_n)}.\]
Then the result follows from Lemma \ref{lem:hgconv}.
\end{proof}

\begin{lemma}
Let $\epsilon$ be such that $0<\epsilon<\frac{1}{4}$.
Under the conditions (A0) and (A1), for all $x$ and $n$
\begin{align}\label{eq:half}
\frac{v_n}{g(A_n)}\le 1\Rightarrow-\epsilon\frac{|x|}{h^{-1}\circ g(A_n)}+\epsilon\frac{h(x)-v_n}{g(A_n)}\ge-\frac{1}{2}.
\end{align}
\end{lemma}

\begin{proof}
By (A1), for $x\ge h^{-1}\circ g(A_n)$,
\[\frac{h(x)}{x}\ge\frac{g(A_n)}{h^{-1}\circ g(A_n)}.\]
It follows that
\[\frac{h(x)}{g(A_n)}-\frac{x}{h^{-1}\circ g(A_n)}\ge0.\]
For $0\le x<h^{-1}\circ g(A_n)$,
\[\frac{h(x)}{g(A_n)}-\frac{x}{h^{-1}\circ g(A_n)}>-1.\]
Therefore for all $x\ge0$ we have
\[\frac{h(x)-v_n}{g(A_n)}-\frac{x}{h^{-1}\circ g(A_n)}>-1-\frac{v_n}{g(A_n)}.\]
\end{proof}

\begin{lemma}
Under the conditions (A0)-(A2), in UFGH Skeptic can force
\begin{align*}
\sum_n \frac{v_n}{g(A_n)}<\infty\Rightarrow \sum_n\frac{x_n}{h^{-1}\circ g(A_n)}\textrm{ converges.}
\end{align*}
\end{lemma}

\begin{proof}
We consider the following two strategies simultaneously such that $\cK_0=1$ and
\[M_n=\pm \epsilon \cK_{n-1}\frac{1}{h^{-1}\circ g(A_n)},\ V_n=\epsilon \cK_{n-1}\frac{1}{g(A_n)}\]
if $v_n/g(A_n)\le1$ and $M_n=V_n=0$ otherwise, where
$0 < \epsilon < 1/4$ is fixed. We denote the corresponding two
capital processes as $\cK_n^{\pm}$.
Since $\sum_n v_n/g(A_n)<\infty$, we have $v_n/g(A_n)\le1$ for all but finite $n$.
Hence we can assume that the inequality in the right-side hand of (\ref{eq:half}) holds for all $n$.
Then
\[\cK_n^{\pm}
=\cK_{n-1}^{\pm}(1\pm \epsilon\frac{x_n}{h^{-1}\circ g(A_n)}+\epsilon\frac{h(x_n)-v_n}{g(A_n)})\ge\frac{1}{2}\cK_{n-1}^{\pm}.\]
Hence $\cK_n^{\pm}>0$ for all $n$.
By the game-theoretic martingale convergence theorem, both $\cK_n^{\pm}$ converge to a non-negative finite value.

By the inequality $t\ge\log(1+t)\ge t-t^2$ for all $t\ge-1/2$,
\begin{align*}
&\epsilon\sum_{k=1}^n\left(\frac{\pm x_n}{h^{-1}\circ g(A_k)}+\frac{h(x_k)-v_k}{g(A_k)}\right)\\
\ge&\log \cK_n^{\pm}\\
\ge&\ \epsilon\sum_{k=1}^n\left(\frac{\pm x_k}{h^{-1}\circ g(A_k)}+\frac{h(x_k)-v_k}{g(A_k)}\right)
-\epsilon^2\sum_{k=1}^n\left(\frac{\pm x_k}{h^{-1}\circ g(A_k)}+\frac{h(x_k)-v_k}{g(A_k)}\right)^2.
\end{align*}
Notice that each of the following infinite sums is finite;
\begin{align*}
\sum_n\frac{h(x_n)}{g(A_n)},\ \sum_n\frac{v_n}{g(A_n)},\ \sum_n\frac{x_n^2}{(h^{-1}\circ g(A_n))^2},\ \sum_n\frac{h(x_n)^2}{g(A_n)^2},\ \sum_n\frac{v_n^2}{g(A_n)^2}.
\end{align*}
By the inequality
\[(a_1+\cdots+a_m)^2\le m(a_1^2+\cdots+a_m^2)\]
we have
\[\sum_{k=1}^n\left(\frac{\pm x_k}{h^{-1}\circ g(A_k)}+\frac{h(x_k)-v_k}{g(A_k)}\right)^2<\infty.\]
It follows that
\[\sup_n\sum_{k=1}^n\frac{\pm x_n}{h^{-1}\circ g(A_k)}<\infty.\]
Hence both $\log \cK_n^{\pm}$ converge to a finite value.
Therefore we obtain the desired result.
\end{proof}

\subsection{The case of divergence}\label{subsec:combdiv}
We now consider the case of divergence of $\sum_n v_n/g(A_n)$.

\begin{theorem}\label{th:combdiv}
Suppose that $h$ satisfies (A0)-(A2) and that $g$ is a positive increasing function.
Then in UFGH Reality strongly complies with 
\begin{align}\label{eq:combdiv}
\sum_n \frac{v_n}{g(A_n)}=\infty\Rightarrow \sum_{k=1}^n\frac{x_k-m_k}{h^{-1}\circ g(A_k)}\textrm{ does not converge.}
\end{align}
\end{theorem}

\begin{proof}
The proof is similar to that of Theorem \ref{th:notconv}.
We will consider the protocol such that Reality's move is restricted for each $n$.
We shall show that Skeptic can force (\ref{eq:combdiv}) and Reality strongly complies 
with (\ref{eq:combdiv})  in this protocol.

Let $y_n=x_n-m_n$, $g_n=h^{-1}\circ g(A_n)$ and $z_n=y_n/g_n$.
We restrict $z_n\in\{0,\pm1\}$.
If $|z_n|=1$ then $h(y_n)=h(g_n)=g(A_n)$.
If $z_n=0$ then $h(y_n)=h(0)=0$.
In any cases we have $h(y_n)= g(A_n)z_n^2$.

The capital process is
\[\cK_n=\cK_{n-1}+M_ng_n z_n+V_ng(A_n)\left(z_n^2-\frac{v_n}{g(A_n)}\right).\]
By Corollary \ref{cor:threepoint}, Skeptic can force
\[
\sum_n\frac{v_n}{g(A_n)}=\infty\Rightarrow |z_n|=1\textrm{ for infinitely many }n
\]
and Reality can choose $x_n$ such that $\cK_n\le \cK_{n-1}$ and hence  $\sup_n \cK_n\le 1$.
If $|z_n|=1$ for infinitely many $n$, then $\sum_n y_n/g_n$ does not converge.
Hence Skeptic can force (\ref{eq:combdiv}).
\end{proof}

We will show that $\frac{\sum_{k=1}^n(x_k-m_k)}{h^{-1}\circ g(A_n)}$ does not converge 
under the additional condition $h(xy)=h(x)h(y)$ for all $x,y\ge 0$.
However this condition implies $h(x)=x^r$ and $h^{-1}(x)=x^{1/r}$ for $x\ge 0$.
Taking into account of (A1)-(A2) we have $1\le r\le 2$.

\begin{theorem}
Let $h(x)=x^r$ where $1\le r\le 2$ and $g$ be a positive increasing function.
Then Reality strongly complies with 
\begin{align*}
\sum_n \frac{v_n}{g(A_n)}=\infty\Rightarrow \frac{\sum_{k=1}^n(x_k-m_k)}{h^{-1}\circ g(A_n)}\textrm{ does not converge.}
\end{align*}
\end{theorem}

\begin{proof}%

Let $y_n=x_n-m_n$, $g_n=h^{-1}\circ g(A_n)$ and $z_n=y_n/g_n$.
This time let $\{\epsilon_n\}$ be such that $h(\epsilon_n)$ is a sequence of Lemma \ref{lem:ser} for $\frac{v_n}{g(A_n)}$.
We can assume that $0 < \epsilon_n\le1$ for all $n$.
We restrict $\epsilon_n z_n\in\{0,\pm1\}$.
Then we have
\begin{align*}
h(y_n)-v_n
=&h(g_n)h(z_n)-v_n
=\frac{g(A_n)}{h(\epsilon_n)}\left(h(\epsilon_n z_n)-h(\epsilon_n)\frac{v_n}{g(A_n)}\right)\\
=&\frac{g(A_n)}{h(\epsilon_n)}\left((\epsilon_n z_n)^2-h(\epsilon_n)\frac{v_n}{g(A_n)}\right).
\end{align*}
Hence the capital process is
\begin{align*}
\cK_n
&=\cK_{n-1}+M_n y_n+V_n(h(y_n)-v_n)\\
&=\cK_{n-1}+M_n \frac{g_n}{\epsilon_n} \epsilon_n z_n
+V_n \frac{g(A_n)}{h(\epsilon_n)}\left((\epsilon_n z_n)^2-h(\epsilon_n)\frac{v_n}{g(A_n)}\right).
\end{align*}
By Corollary \ref{cor:threepoint}, Skeptic can force 
\begin{align*}
\sum_nh(\epsilon_n)\frac{v_n}{g(A_n)}=\infty \Rightarrow \epsilon_n |z_n|=1\textrm{ for infinitely many }n
\end{align*}
and Reality can choose $x_n$ such that $\cK_n\le \cK_{n-1}$ and hence $\cK_n \le 1$.
The rest of the proof is the same as that of Theorem \ref{th:notconv}.
\end{proof}

One may think that this result contradicts Marcinkiewicz-Zygmund strong law,
which says that for i.i.d.\ random variables $\{x_n\}$ with $E|x_n|^r<\infty$ for $0<r<2$ and $Ex_n=0$ when $1\le r<2$,
$n^{-1/r}(\sum_{k=1}^n x_k)\to0$ as $n\to\infty$ a.s.

For example let $m_n=0$, $v_n=v$ for all $n$ and $g(x)=x/v$.
Then $g_n=(g(A_n))^{1/r}=n^{1/r}$ and $\epsilon_n z_n=\epsilon_n y_n/g_n=\epsilon_n x_n/n^{1/r}$.
By $\epsilon_n z_n\in\{0,\pm 1\}$ we have $x_n\in\{0,\pm n^{1/r}/\epsilon_n\}$.
Since $\epsilon_n\to0$, we have $n^{1/r}/\epsilon_n\to\infty$ as $n\to\infty$.
It follows that the restrictions of $x_n$ are not the same.
Thus this is not a contradiction.

Theorem \ref{th:comb} and Theorem \ref{th:combdiv} give the rate of 
convergence of SLLN under a general hedge in a non-identical case.

\section{Discussion}\label{sec:disc}

In the classical
three-series theorem, under the assumption of 
independence of random variables, the necessary and sufficient condition
for the convergence of a random series is given by the
convergence of three series: truncation probabilities,
truncated expected values and truncated variances.
On the other-hand Gilat's counter example shows
that the necessity of the convergence of three series can
not hold for martingales.  Therefore a question of interest is
to specify some conditions, other than the independence, under which
the convergence of a random series implies the convergence of three series.
We have shown that under the convergence of a random series, 
the divergence of truncated means can only occur
as two-sided unbounded oscillation. We also gave 
some simple conditions of convergences of truncated means and variances,
but stronger results are desirable.

We proposed the notion of compliance concerning Reality's deterministic 
strategy.  We showed that a good deterministic strategy of Reality can be 
automatically constructed
by using a good strategy for Skeptic as a ``surrogate''.  
In Definition \ref{def:reality-stragety} we made a distinction
between compliance and strong compliance concerning
Reality's deterministic strategy.  Further study is needed
to clarify the difference between these definitions.

We gave the precise limit of the rate of convergence of SLLN with the quadratic hedge as well as more general weaker hedges.
According to it the rate of convergence for a non-identical case may be slower than for the i.i.d.\ case.
We believe that game-theoretic probability is a powerful tool for analysis of such a non-identical case.

\medskip

\noindent {\bf Acknowledgment}:

We thank a reviewer for very detailed and constructive comments.
The first author is supported by GCOE, Kyoto University.
The second author is partially supported by the Aihara Project, the FIRST
program from JSPS.
We are grateful to Hideatsu Tsukahara for useful comments.

\bibliographystyle{abbrv}
\bibliography{threeseries}

\begin{thebibliography}{10}

\bibitem{AroBar09}
S.~Arora and B.~Barak.
\newblock {\em {Computational Complexity: a Modern Approach}}.
\newblock Cambridge University Press, 2009.

\bibitem{brown1971}
B.~M. Brown.
\newblock A general three-series theorem.
\newblock {\em Proc. Amer. Math. Soc.}, 28:573--577, 1971.

\bibitem{brown1972erratum}
B.~M. Brown.
\newblock Erratum: ``{A} general three-series theorem''.
\newblock {\em Proc. Amer. Math. Soc.}, 32:634, 1972.

\bibitem{Dow10}
R.~Downey and D.~R. Hirschfeldt.
\newblock {\em Algorithmic Randomness and Complexity}.
\newblock Springer, Berlin, 2010.

\bibitem{fisher-1992}
E.~Fisher.
\newblock On the law of the iterated logarithm for martingales.
\newblock {\em Ann. Probab.}, 20(2):675--680, 1992.

\bibitem{Gac85}
P.~G\'acs.
\newblock Every set is reducible to a random one.
\newblock {\em Information and Control}, 70:186--192, 1986.

\bibitem{gilat1971}
D.~Gilat.
\newblock On the nonexistence of a three series condition for series of
  nonindependent random variables.
\newblock {\em The Annals of Mathematical Statistics}, 42(1):409, 1971.

\bibitem{gut2005}
A.~Gut.
\newblock {\em Probability: a Graduate Course}.
\newblock Springer, New York, 2005.

\bibitem{HarWin41}
P.~Hartman and A.~Wintner.
\newblock On the law of the iterated logarithm.
\newblock {\em American J. Math.}, 63:169--176, 1941.

\bibitem{horikoshi-takemura-2008}
Y.~Horikoshi and A.~Takemura.
\newblock Implications of contrarian and one-sided strategies for the fair-coin
  game.
\newblock {\em Stochastic Process. Appl.}, 118(11):2125--2142, 2008.

\bibitem{Khi24}
A.~Y. Khinchin.
\newblock {\"U}ber einen {S}atz der {W}ahrscheinlichkeitrechnung.
\newblock {\em Fund. Mat.}, 6:9--20, 1924.

\bibitem{Kol29}
A.~N. Kolmogorov.
\newblock {\"Uber} das {G}esetz des {I}terierten {L}ogarithmus.
\newblock {\em Math. Ann.}, 101:126--135, 1929.

\bibitem{KumTakTak07}
M.~Kumon, A.~Takemura, and K.~Takeuchi.
\newblock {Game-theoretic versions of strong law of large numbers for unbounded
  variables}.
\newblock {\em Stochastics}, 79(5):449--468, 2007.

\bibitem{KumTakTak11}
M.~Kumon, A.~Takemura, and K.~Takeuchi.
\newblock {Sequential optimizing strategy in multi-dimensional bounded
  forecasting games}.
\newblock {\em Stochastic Processes and their Applications}, 121:155--183,
  2011.

\bibitem{marcinkeiwicz-zygmund}
J.~Marcinkiewicz and A.~Zygmund.
\newblock Sur les fonctions ind\'ependantes.
\newblock {\em Fund. Math.}, 29:60--90, 1937.

\bibitem{MerMih02}
W.~Merkle and N.~Mihailovi{\'c}.
\newblock {On the construction of effective random sets}.
\newblock {\em Mathematical Foundations of Computer Science}, pages 568--580,
  2002.

\bibitem{Miyabe_extvan}
K.~Miyabe.
\newblock An extension of van {L}ambalgen's {T}heorem to infinitely many
  relative 1-random reals.
\newblock {\em Notre Dame Journal of Formal Logic}, 51(3):337--349, 2010.

\bibitem{nakajima-etal}
R.~Nakajima, M.~Kumon, A.~Takemura, and K.~Takeuchi.
\newblock Approximations and asymptotics of upper hedging prices in multinomial
  models, 2010.
\newblock {\tt arXiv:1007.4372v1}. To appear in {\it Japan Journal of
  Industrial and Applied Mathematics}.

\bibitem{Nie09}
A.~Nies.
\newblock {\em Computability and Randomness}.
\newblock Oxford University Press, USA, 2009.

\bibitem{NST05}
A.~Nies, F.~Stephan, and S.~Terwijn.
\newblock Randomness, relativization and {T}uring degrees.
\newblock {\em Journal of Symbolic Logic}, 70:515--535, 2005.

\bibitem{Odi90}
P.~Odifreddi.
\newblock {\em Classical Recursion Theory}, volume~1.
\newblock North-Holland, 1990.

\bibitem{Odi99}
P.~Odifreddi.
\newblock {\em Classical Recursion Theory}, volume~2.
\newblock North-Holland, 1999.

\bibitem{Pet95}
V.~V. Petrov.
\newblock {\em Limit Theorems of Probability Theory: Sequences of Independent
  Random Variables}.
\newblock Oxford University Press, USA, 1995.

\bibitem{ShaVov01}
G.~Shafer and V.~Vovk.
\newblock {\em Probability and Finance: It's Only a Game!}
\newblock Wiley, 2001.

\bibitem{takazawa-aism}
{Shin-ichiro Takazawa}.
\newblock Exponential inequalities and the law of the iterated logarithm in the
  unbounded forecasting game, 2010.
\newblock To appear in {\it Annals of the Institute of Statistical
  Mathematics}.

\bibitem{takazawa-stochastics}
{Shin-ichiro Takazawa}.
\newblock An exponential inequality and the convergence rate of the strong law
  of large numbers in the unbounded forecasting game.
\newblock {\em Stochastics}, 83:117--125, 2011.

\bibitem{Shi95}
A.~N. Shiryaev.
\newblock {\em {Probability}}.
\newblock Springer, second edition, 1995.

\bibitem{stout-1970}
W.~F. Stout.
\newblock A martingale analogue of {K}olmogorov's law of the iterated
  logarithm.
\newblock {\em Z. Wahrscheinlichkeitstheorie und Verw. Gebiete}, 15:279--290,
  1970.

\bibitem{tomkins-1978}
R.~J. Tomkins.
\newblock On the law of the iterated logarithm.
\newblock {\em Ann. Probability}, 6(1):162--168, 1978.

\bibitem{vovk-shen-2010}
V.~Vovk and A.~Shen.
\newblock Prequential randomness and probability.
\newblock {\em Theoret. Comput. Sci.}, 411(29-30):2632--2646, 2010.

\bibitem{Wil91}
D.~Williams.
\newblock {\em Probability with Martingales}.
\newblock Cambridge University Press, 1991.

\end{thebibliography}

\end{document}